\documentclass[final,leqno]{siamltex}

\usepackage[dvips]{color}

\usepackage{cite}
\usepackage{tabularx}

\usepackage{pgf}

\usepackage{amsmath,amssymb}
\usepackage{amsmath}
\usepackage{verbatim}
\usepackage{graphicx}
\graphicspath{{../figs/}}
\usepackage{algorithm}
\usepackage{algorithmic}
\usepackage{amsmath,amscd}
\usepackage[all]{xy}
\usepackage{float}

\newcommand{\bbm}{\begin{bmatrix}}
\newcommand{\ebm}{\end{bmatrix}}
\newcommand{\beq}{\begin{equation}}
\newcommand{\beqn}{\begin{equation*}}
\newcommand{\eeq}{\end{equation}}
\newcommand{\eeqn}{\end{equation*}}

\newtheorem{remark}[theorem]{Remark}
\numberwithin{equation}{section}

\newcommand{\tA[1]}{\widetilde{#1}}

\newcommand{\cA[1]}{\mathcal{#1}}
\newcommand{\ttA[1]}{\mathtt{#1}}
\newcommand{\acc}{\varepsilon}

\title{An O(N) Direct Solver for Integral Equations on the Plane}
\author{Eduardo Corona \footnotemark[2] \and Per-Gunnar Martinsson \footnotemark[3] 
\and Denis Zorin \footnotemark[2]}

\begin{document}
\maketitle

\renewcommand{\thefootnote}{\fnsymbol{footnote}}
\footnotetext[2]{Courant Institute of Mathematical Sciences, New York University}
\footnotetext[3]{Department of Applied Mathematics, University of Colorado at Boulder}

\bibliographystyle{plain}

\begin{abstract}
An efficient direct solver for volume integral equations
with $O(N)$ complexity for a broad range of problems is presented.
The solver relies on hierarchical compression of the discretized integral operator,
and exploits that off-diagonal blocks of certain dense matrices have numerically low rank.
Technically, the solver is inspired by previously developed direct solvers for integral
equations based on ``recursive skeletonization'' and ``Hierarchically Semi-Separable'' (HSS)
matrices, but it improves on the asymptotic complexity of existing solvers by
incorporating an additional level of compression. The resulting solver has optimal
$O(N)$ complexity for all stages of the computation, as demonstrated
by both theoretical analysis and numerical examples. The computational examples further
display good practical performance in terms of both speed and memory usage. In particular, it is demonstrated
that even problems involving $10^{7}$ unknowns can be solved to precision $10^{-10}$ using
a simple Matlab implementation of the algorithm executed on a single core.
\end{abstract}

\section{Introduction}
\label{sec:introduction}
Many boundary value problems from classical physics, when cast as boundary or volume integral equations,
take the form
\beq
a(x) \sigma(x) + \int_{\Gamma} {\cA[K](x,y) \sigma(y) dS(y)} = f(x),\qquad  \forall x \in \Gamma,
\label{eq:basic-integral}
\eeq
where $\Gamma$ is a domain in either $\mathbb{R}^{2}$ or $\mathbb{R}^{3}$ (either
a boundary or a volume) and $\cA[K](x,y)$ is a \textit{kernel function}
derived from the fundamental solution
associated with the relevant elliptic PDE (e.g.~the Laplace or Helmholtz equation,
the Stokes equations, etc.). The kernel function is typically singular near
the diagonal (as $y$ approaches $x$) but is otherwise smooth.

Discretizing the integral in (\ref{eq:basic-integral}) using, e.g., the
Nystr\"om method \cite{atkinson1997} results in a linear system
of the form
\beq
A\sigma = f,
\label{eq:Asigma}
\eeq
where $A$ is a dense $N \times N$ matrix.
For many problems it has been known since the 1980's how to rapidly
evaluate the matrix-vector product $\sigma \mapsto A\sigma$ using $O(N)$
operations, see, e.g., \cite{greengard1987fast}. Such matrix-vector multiplication
techniques can be coupled with an iterative solver
(e.g.~GMRES~\cite{GMRES} or Bi-CGSTAB~\cite{BiCGSTAB})
to attain solvers for (\ref{eq:basic-integral}) that can be
very effective whenever convergence is rapid. More recently,
\textit{direct solvers} with linear complexity have been developed
for (\ref{eq:Asigma}) for the case when $\Gamma$ is a contour in the
plane \cite{MR2005}. When the method of \cite{MR2005} is applied
to a volume domain in the plane, a reasonably efficient direct solver
with asymptotic complexity $O(N^{3/2})$ results \cite{GYMR2012,ho2012fast,G2011thesis}.
This paper describes how the algorithm of \cite{MR2005} can be
reworked to construct a direct solver for volume integral equations
in two dimensions that has optimal $O(N)$ complexity and high
practical efficiency (even at high accuracies such as $10^{-10}$).
We believe that this solver will also (with some additional work)
be capable of directly solving boundary integral equations for
problems in 3D in $O(N)$ or $O(N\log N)$ operations, and will be
very well suited for implementation on multicore and parallel
computers.

\subsection{Previous Work}

The key observation enabling the construction of $O(N)$ algorithms for
solving linear systems arising from the discretization of integral
equations is that the off-diagonal blocks of the coefficient matrix
can be very well approximated by matrices of low numerical rank. In
this section, we briefly describe some key results.

\subsubsection*{Optimal complexity iterative solvers for integral equations}
By combining the observation that off-diagonal blocks of the matrix
have low rank with a hierarchical partitioning of the physical domain
into a tree-structure, algorithms of $O(N)$ or $O(N\log N)$ complexity
for matrix-vector multiplication were developed in the 1980's.
Perhaps, the most prominent is the Fast Multiple Method~\cite{greengard1987fast,rokhlin1985,rokhlin1997},
but the Panel Clustering ~\cite{PanelClustering} and Barnes-Hut~\cite{BarnesHut}
methods are also well known. When a fast algorithm for the matrix-vector
multiplication such as the FMM is coupled with an iterative method such as, e.g.,
GMRES~\cite{GMRES} or Bi-CGSTAB~\cite{BiCGSTAB}, the result is a solver
for the linear system arising upon discretization of (\ref{eq:basic-integral})
with overall complexity $O(M\,N)$, where $M$ is the number of steps
required by the iterative solver. In many situations, $M$ is independent
of $N$ and convergence can be very rapid.

The first FMMs constructed were custom designed for specific elliptic
equations (e.g.~Laplace, Helmholtz, Stokes), but it was later realized
that \textit{kernel-independent methods} that work for broad classes of
problems could be developed~\cite{kifmm04ying,gimbutas2002}.
 They directly inspired the work described in this paper, and
the new direct solvers can also be said to be ``kernel-independent'' in the
sense that the same algorithm, and the same code, can be applied to several
different types of physical problems such as electro-statics, Stokes
flows, and low frequency scattering.

While iterative solvers accelerated by fast methods for matrix-vector
multiply can be very effective in many contexts, their performance
is held hostage to the convergence rate of the iteration. If the equation is not well-conditioned, the
complexity of an iterative solve may increase. Pre-conditioners can sometimes be constructed to accelerate convergence,
but these tend to be quite problem-specific and do not readily lend themselves
to the construction of general purpose codes.
Examples of ill-conditioned problems that can be challenging to iterative methods include
Fredholm equations of the first kind, elasticity problems on thin domains,
and scattering problems near resonances.

\subsubsection*{Direct solvers for integral equations}
In the last ten years, a number of efficient \textit{direct} solvers for
linear systems associated with integral equations have been constructed.
These solvers entirely side-step the challenges related to convergence
speed of iterative solvers. They can also lead to dramatic improvements
in speed, in particular in situations where several of equations with
the same coefficient matrix but different right-hand sides need to be solved.

The work presented here draws heavily on \cite{MR2005}, (based on
on \cite{starr_rokhlin}), which describes a direct solver that was
originally developed for boundary integral equations defined on curves
in the plane, and has optimal $O(N)$ complexity for this case. The
observation that this algorithm
can also be applied to volume integral equations in the plane, and to
boundary integral equations in 3D was made in \cite{greengard2009fast}
and later elaborated in \cite{G2011thesis}.
For these cases, the direct solver requires $O(N^{3/2})$ flops to build
an approximation to the inverse of the matrix, and $O(N\log N )$ flops for the
``solve stage'' once the inverse has been constructed.
Similar work was done in \cite{ho2012fast}, where it is also demonstrated that
the direct solver can, from a practical point of view, be implemented using
standard direct solvers for large \textit{sparse matrices}. This improves
stability, and greatly simplifies the practical implementation due to the
availability of standard packages such as UMFPACK.

The direct solver of \cite{MR2005} relies on the fact that the matrices
arising form the discretization of integral equations can be efficiently
represented in a data-sparse format often referred to as ``Hierarchically
Semi-Separable (HSS)'' matrices. This matrix format was also explored in
\cite{Gu06sparse,Gu06ULV}, with a more recent efficient version presented in
\cite{xia2010}. This work computes ULV and Cholesky factorizations of HSS
matrices; if these techniques were applied to volume integral equations in
2D, the complexity would be $O(N^{3/2})$, in complete agreement with
\cite{greengard2009fast} and \cite{G2011thesis}.
The paper \cite{xia2012complexity} presents a general complexity study of
HSS algorithms, under different rank growth patterns; it presents an
optimal-complexity  HSS recompression method which we adapt to our
setting as a part of the overall algorithm.

An important class of related algorithms is  $\cA[H]$ and $\cA[H]^2$-matrix
methods of Hackbusch and co-workers (see \cite{borm2003hierarchical,2008_bebendorf_book,2010_borm_book}
for surveys). These techniques are based on variations of the cross approximation method for low-rank compression,
and have been applied both to integral equations and sparse systems derived from PDEs.
The matrix factorization for two and three-dimensional problems algorithms are formulated
recursively, and a full set of compressed operations for lower-dimensional problems needs to be
available. In \cite{borm2009construction,borm2006matrix}, algorithms for
$\cA[H]^2$ matrix arithmetics are described;  Observed behavior for integral equation operators on the cube and on the sphere in chapter 10 of \cite{2010_borm_book} is $O(N\log^4 N )$ for matrix compression, $O(N\log^3 N )$ for inversion and $O(N\log^2 N )$ for solve time and memory use.

\subsubsection*{Direct solvers for sparse systems}
Our direct solver  is conceptually  related to direct solvers
for sparse system matrices such as the classical \textit{nested dissection} and
\textit{multifrontal} methods \cite{george_1973,hoffman_1973,1989_directbook_duff}.
These solvers do not have optimal complexity (they typically require $O(N^{3/2})$
for the factorization stage in 2D, and $O(N^{2})$ in 3D), but are nevertheless
popular (especially in 2D) due in part to their robustness, and in part to the
unrivaled speed that can be attained for problems involving multiple right hand
sides. Very recently, it has been demonstrated that by exploiting structured matrix
algebra (such as, e.g., $\mathcal{H}$-matrices, or HSS matrices), to manipulate the
dense matrices that arise due to fill-in, direct solvers of linear or close to linear
complexity can be constructed
\cite{2009_xia_superfast,2011_ying_nested_dissection_2D,2007_leborne_HLU,2009_martinsson_FEM,G2011thesis}.
The direct solver described in this paper is conceptually similar to these algorithms
in that they all rely on hierarchical domain decompositions, and efficient representations
of operators that live on the interfaces between sub-domains.

\subsection{Overview of new results}

We present a direct solver that achieves optimal $O(N)$ complexity for \emph{two-dimensional} systems derived from integral equations with non-oscillatory kernels or kernels in low-frequency mode. Similarly to other HSS or $\cA[H]$-matrix methods, we rely on the fact that the system is derived from an integral equation only weakly: if the kernel is of a different nature, the scalability of the solver may deteriorate, but it can still perform accurate calculations.

The main features of our solver include:

\begin{itemize}
\item Observed $O(N)$ complexity both in time and storage for all stages of the computation
(in particular, for both the ``build'' and the ``solve'' stages).
\item The algorithm supports high accuracy (up to $10^{-10} - 10^{-12}$ is practical), while maintaing reasonable efficiency both in time and memory.
\item The algorithm can take direct advantage of translation invariance and symmetry of underlying kernels, achieving considerable speedup and reduction of memory cost.
\end{itemize}

The main aspects of the algorithm that allow us to achieve high performance are:

\begin{itemize}
\item Two levels of hierarchical structures are used. The matrix $A$ as whole is represented in the HSS format,
and  certain blocks within the HSS structure are themselves represented in the HSS format.
\item Direct construction of the inverse: unlike many previous algorithms, we directly build a compressed representation of the inverse, rather than compressing the matrix itself and then inverting.
\item Our direct solver needs only a subset of a full set of HSS matrix arithmetics;
in particular, the relatively expensive matrix-matrix multiplication is never used.
\end{itemize}

We achieve significant gains in speed and memory efficiency with
 respect to the existent $O(N^{3/2})$ approach. For non-oscillatory
 kernels, our algorithm outperforms the $O(N^{3/2})$ algorithm around
 $N \sim10^5$. Sizes of up to $N \sim 10^7$ are practical on a desktop
 gaining one order of magnitude in inverse compression time and
 storage. For example, for non-translation-invariant kernels, a 
 problem of size $N = 3 \times 10^6$ and target accuracy $\acc = 10^{-10}$ 
 takes 1 hour and $\sim$ 50GB to invert. Each solve takes $10$ seconds.
For translation-invariant kernels such as Laplace, a problem of size
$N = 1 \times 10^7$ and target accuracy $\acc = 10^{-10}$ takes half
an hour and $5$ GB to invert, with $20$-second solves.
Reducing target accuracy to $\acc = 10^{-5}$, inversion costs for the 
latter problem go down to 5 minutes and $1$GB of storage, and each solve takes $10$ seconds.

Our accelerated approach to build the HSS binary tree also yields a fast $O(N)$ matrix compression algorithm. As is noted in \cite{ho2012fast}, given that matrix-vector applies are orders of magnitude faster than one round of FMM, this algorithm would be preferable to an interative method coupled with FMM for problems that require more than a few iterations. To give an example, for $N=10^7$ and $\acc = 10^{-10}$, matrix compression takes $5$ minutes and $1$GB of storage, and each matrix apply takes less than $10$ seconds.

Finally, we are able to apply our method with some minor modifications to oscillatory kernels in low frequency mode, and apply this to the solution of the corresponding 2D volume scattering problems. Although the costs are considerably higher, we still observe optimal scaling and similar performance gains.

\section{Background}
\label{sec:background}

Our approach builds on the fast direct solver in \cite{GYMR2012}, which consists of separate
hierarchical compression, inversion and inverse apply algorithms, all of which achieve linear
complexity when applied to integral operators on one dimensional curves.
These algorithms can be applied to volume integral equations in 2D and to
boundary integral equations defined on surfaces in 3D but compression and
inversion then have $O(N^{3/2})$ complexity, while application of the inverse is $O(N \log N)$.
The super-linear complexity renders the technique impractical for large
problems, it is in practice useful for small and moderate size problems,
see \cite{ho2012fast,G2011thesis}.

In this section we review the algorithms for compression and inversion described in \cite{GYMR2012}, as these are needed to describe our algorithms, and discuss their limitations for problems on 2D domains.

\subsection{Notation and Preliminaries}

We view an $N \times N$ matrix $A$  as a kernel function $K = K(p,q)$ evaluated at pairs of sample points. As it typically comes from an integral formulation of an elliptic PDE, aside from a low-rank block structure we also expect Green's identities to hold for $K$; this however, is not a fundamental limitation of our algorithm: the use of Green's identities is restricted to one small part of the algorithm (Section~\ref{sec:hss}) relying on equivalent density representation and can be replaced by any other technique of similar nature.

We use matlab-like notation $A(I,J)$ where $I$ and $J$ are ordered sets of
indices to denote submatrices of a matrix $A$. Again, following
matlab conventions, $A(:,J)$ and $A(I,:)$ indicate blocks of columns
and rows, respectively.

An essential building block for our algorithm is the interpolative
low-rank decomposition. This decomposition factors an $m \times n$
matrix $A$ into a narrower \emph{skeleton} matrix $A^{\rm sk} = A(:, I^{\rm sk})$
of size $m \times k$,  consisting of a subset of columns of $A$
indexed by the set of indices  $I^{\rm sk}$, and the interpolation matrix
$T$ of size $k \times (n-k)$, expressing the remaining columns of $A$
as linear combinations of columns of $A^{\rm sk}$: The set of indices
$I^{\rm sk}$ is called the \emph{column skeleton} of $A$. If $\Pi^{\rm sk}$ is
the permutation matrix placing entries with indices from $I^{\rm sk}$ first,

\begin{equation}
A = A^{\rm sk} \bbm I_{k \times k}\; T \ebm \Pi^{\rm sk} + E
\label{eq:interp-fact-R}
\end{equation}
where $||E||_2 \sim \sigma_{k+1}$ vanishes as we increase $k$, and $R
= [I_{k\times k} T] \Pi^{\rm sk}$ is a downsampling interpolation matrix.
We denote this compression operation by $[T,I^{\rm sk}] = ID(A, \acc)$, where
$\acc$ controls the norm of $E$. To obtain a similar compression for
rows, we apply the same operation to $A^T$; in this case, we obtain a
factorization $A =  (\Pi^{\rm sk}_{r})^T \bbm I_{k \times k} \\ T_r^T \ebm A^{\rm sk}_r + E_r$,
where $L = (\Pi^{\rm sk}_r)^T \bbm I_{k \times k} \\ T_r^T \ebm$
is an upsampling interpolation matrix.

\subsection{Constructing hierarchically semi-separable matrices}
\label{sec:hss}

We assume that the domain of interest is contained in a  rectangle
$\Omega$, with a regular grid of samples (if the input data is given
in a different representation, we resample it first).
A quadtree is constructed by recursively subdividing  $\Omega$ into
cells (boxes $B_i$)  by bisection,
corresponding to  the nodes of a \emph{binary} tree $\cal T$.
The two \emph{children} of a box are denoted $c_1(i)$ and $c_2(i)$.
Subdivision can be done adaptively without significant changes to the
algorithm,  but to simplify the exposition we focus on the uniform refinement case.
Let $I_i$ denote the index vector marking all discretization points in
box $B_{i}$, and let $\mathcal{L}$ denote a list of the leaf boxes. Then
$\{I_{i}\}_{i \in \mathcal{L}}$ forms a disjoint partition of the full
index set
\begin{equation}
\label{eq:indexpartition}
\{1,\,2,\,3,\,\dots,\,N\} = \bigcup_{i \in \mathcal{L}} I_{i}.
\end{equation}
Let $m_{i}$ denote the number of points in $B_{i}$.
The partition (\ref{eq:indexpartition}) corresponds to a blocking
\begin{equation}
\label{eq:blocked}
\sum_{j \in \mathcal{L}} {A_{ij} \sigma_j} = f_i,\qquad i \in \mathcal{L},
\end{equation}
of the linear system (\ref{eq:Asigma}), where
$A_{ij} = A(I_i,I_j)$, and the vectors $\sigma$ and $f$ are partitioned accordingly.
For a linear system such as (\ref{eq:Asigma}) arising from the discretization
of an integral equation with a smooth kernel, the off-diagonal blocks of
(\ref{eq:blocked}) typically have low numerical rank. Such matrices can be
represented in an efficient data-sparse format called
\emph{hierarchically semi-separable (HSS)}. In order to
rigorously describe this format, we first define the concept of
a \emph{block-separable matrix}.

\begin{definition} [Block-separable Matrices]
\label{def:blockseparable}
We say $A$ is \emph{block-separable} if there exist matrices $\{L_i,R_i \}_{i \in \mathcal{L}}$
such that each off-diagonal block $A_{i,j}$ in (\ref{eq:blocked}) admits the factorization
\begin{equation}
\label{eq:blockseparable}
A_{i,j} = \underset{m_i \times k_i}{L_i} \ \underset{k_i \times k_j}{M_{i,j}} \ \underset{k_j \times m_j}{R_j},
\end{equation}
where the block ranks $k_{i}$ satisfy $k_{i} < m_{i}$.
\end{definition}

In order to construct the matrices $L_{i}$ and $R_{i}$ in (\ref{eq:blockseparable})
it is helpful to introduce  \textit{block rows} and \textit{block columns} of $A$:
For $i \in \mathcal{L}$, we define the $i$th off-diagonal block row of
$A$ as $A^{row}_i = A(I_i, I \setminus I_i) = [A_{i,1} \dots A_{i,i-1}
A_{i,i+1} \dots A_{i,p}]$. The $j$th off-diagonal block column
$A^{col}_j$ is defined analogously. Given a prescribed accuracy $\acc$, we denote by
$k_i^{r}$ and $k_i^{c}$ the $\acc$-ranks of $A^{row}_i$ and $A^{col}_j$, respectively.

Now that we have defined $A^{row}_i$ and $A^{col}_i$, we use them to obtain
the factorization (\ref{eq:blockseparable}) as follows:
For each $i \in \mathcal{L}$, form interpolative decompositions of
$A^{row}_i$ and $A^{col}_i$:
\begin{equation}
\label{eq:fullID}
A_i^{row} = \underset{m_i \times k_i}{L_{i}}
\underset{k_{i} \times (N-m_i)}{A(I_{i}^{\rm rsk},:)}
\qquad\mbox{and}\qquad
A_i^{col} =
\underset{(N-m_i) \times k_i}{A(I_{i}^{\rm csk},:)}
\underset{k_i \times m_i}{R_{i}},
\end{equation}
where the index vectors $I_{i}^{\rm rsk}$ and $I_{i}^{\rm csk}$ are
the \textit{row-skeleton} and \textit{column-skeleton} of block $i$, respectively.
Note that the columns of $L_i$ are a column basis for $A^{row}_i$,
and the rows of $R_i$ a row basis of  $A^{col}_i$. Setting
$$
M_{i,j} = A(I_i^{\rm rsk},I_j^{\rm csk}),
$$
we then find that (\ref{eq:blockseparable}) necessarily holds.
Observe that each matrix $M_{i,j}$ is a submatrix of $A$.

Set $D_i = A_{i,i}$. This yields a block factorization for $A$:
\begin{equation}
A = D^{d} + L^{d} A^{d-1} R^{d}
\label{eq:onelevel-factorization}
\end{equation}
where $D^{d}$, $L^{d}$, and $R^{d}$ are the block diagonal matrices whose diagonal
blocks are given by $\{D_{i}\}_{i \in \mathcal{L}}$, $\{L_{i}\}_{i \in \mathcal{L}}$,
and $\{R_{i}\}_{i \in \mathcal{L}}$, respectively. The matrix $A^{d-1}$ is
the submatrix of $A$ corresponding to the union of skeleton points, with
diagonal blocks zeroed out and off-diagonal blocks $M_{i,j}$.

\begin{remark}
Row and column skeleton sets need not coincide, although for all purposes, we will assume they are augmented so that they are the same size (so $D_i$ blocks are square). If the system matrix is symmetric, as is the case for all matrices considered in this paper, these sets are indeed the same and further $R_i = L_i^{T}$. For this reason, as well as for the sake of simplicity, we make no further distinction between them unless it is necessary.  
\end{remark}

\subsubsection*{Hierarchical compression of A}
\begin{figure}[H]
\begin{center}
\includegraphics[scale = 0.5]{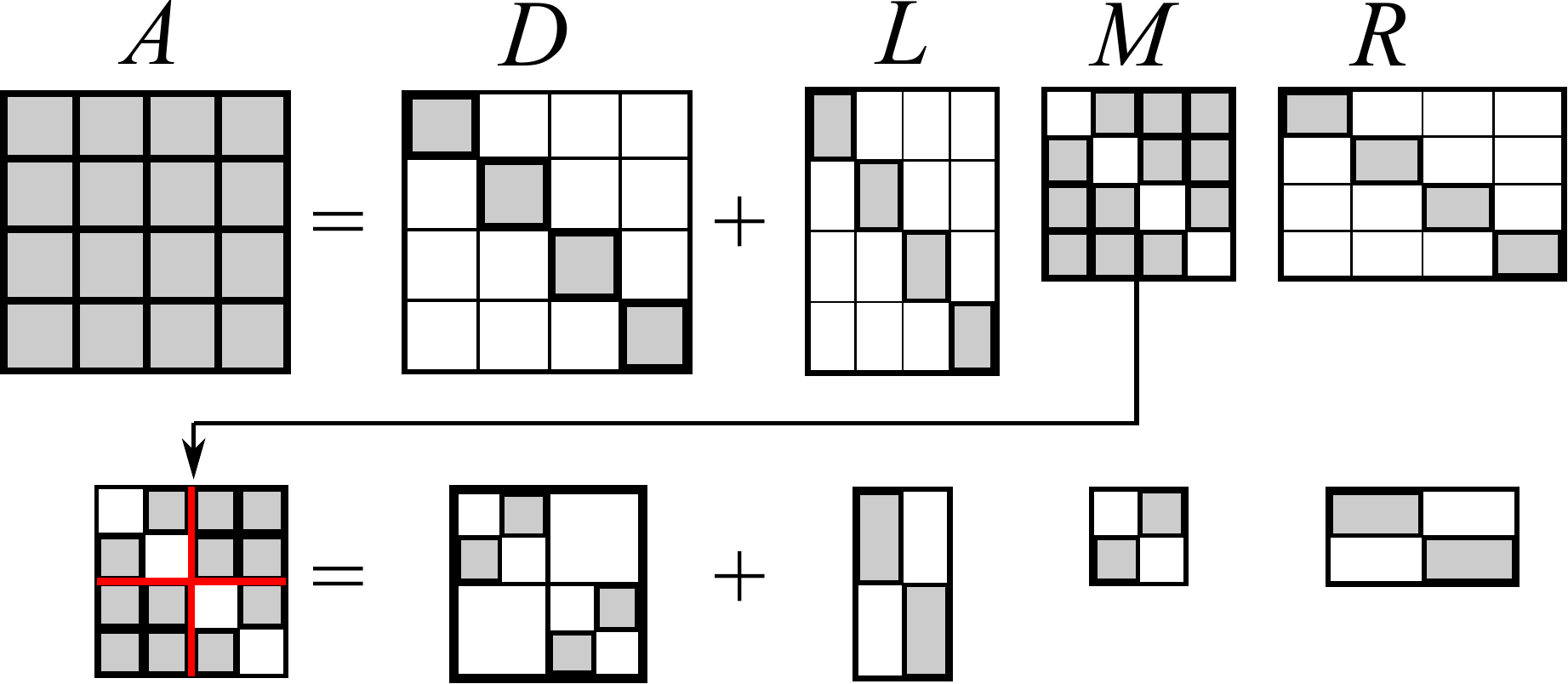}
\caption{\textit{Two levels of block-separable compression:
blocks of $M$ corresponding to children are merged and then
off-diagonal interactions are further compressed.}}
\label{fig:reblock}
\end{center}
\end{figure}

The key property that allows dense matrix operations to be performed
with less that $O(N^{2})$ complexity is that the low-rank structure in
definition~\ref{def:blockseparable} can be exploited recursively in the
sense that the matrix $A^{d-1}$ in \eqref{eq:onelevel-factorization}
itself is block-separable.
To be precise, we re-partition the matrix $A^{d-1}$ by merging $2\times 2$
sets of blocks to form new larger blocks. Each larger block is associated
with a box $i$ on level $d-1$ corresponding to the index vector
$I_{c_1(i)}^{\rm sk} \sqcup I_{c_2(i)}^{\rm sk}$ (the new index vector
holds $m_i = k_{c_1}+k_{c_2}$ nodes).
The resulting matrix with larger blocks is then itself block-separable
and admits a factorization, cf.~Figure \ref{fig:reblock},
\begin{equation}
A = D^{d} + L^{d} \big(D^{d-1} + L^{d-1}\,A^{d-2}\,R^{d-1}\bigr) R^{d}
\label{eq:twolevel-factorization}
\end{equation}

We say that $A$ is \emph{hierarchically semiseparable} (HSS) if the
process of reblocking and recompression can be continued through all
levels of the tree. In other words, we assume that
$A^{\ell} = D^{\ell} + L^{\ell} A^{\ell-1} R^{\ell}$
for $\ell = d \ldots 1$, or, more explicitly,
\begin{equation}
A^d = D^{d} + L^{d}( D^{d-1} + L^{d-1} ( D^{d-1} + \ldots  (D^1 + L^1D^0R^1) \ldots ) R^{d-1}) R^{d}
\label{eq:telescope}
\end{equation}
with $A = A^d$, $A^0 = D^0$, and $D^{\ell}$, $L^{\ell}$ and $R^{\ell}$
block-diagonal, with blocks in matrices with index $\ell$ corresponding
to boxes at level $\ell$. For non-leaf boxes, blocks $D_i^{\ell}$ account
for ``sibling interactions,'' in other words interactions between the
children $c_{1}(i)$ and $c_{2}(i)$ of $B_{i}$, 
\begin{equation}
D_i =
 \left[ \begin{array}{cc}
     & M_{c_1 (i), c_2 (i)}\\
     M_{c_2 (i), c_1 (i)} &
   \end{array} \right] \\
\label{eq:D-expressions}
\end{equation}
We call (\ref{eq:telescope}) \emph{the telescoping factorization} of $A$.
The matrices under consideration in this manuscript are (like most matrices arising
from the discretization of integral operators) all HSS.

\begin{remark}
In this paper we use the term \textit{hierarchically semi-separable (HSS)}
to conform with standard use in the literature, see \cite{Gu06sparse,Gu06ULV,chandrasekaran2010numerical,xia2010}.
In \cite{MR2005,GYMR2012,ho2012fast} the term \textit{hierarchically block-separable (HBS)} is
alternatively used to refer to this hierarchical version of block-separability, consistent
with definition~\ref{def:blockseparable}.
\end{remark}

\subsubsection*{Using equivalent densities to accelerate compression}
A matrix is block-separable as long as all sub-matrices $A^{\rm row}_{i}$
and $A_{i}^{\rm col}$ are rank-deficient. However, these matrices are large,
so directly computing the IDs in (\ref{eq:fullID}) is expensive ($O(N^{2})$ cost \cite{xia2012complexity}).
In this section, we describe how to exploit the fact that the
matrix to be compressed is associated with an elliptic PDE to
reduce the asymptotic cost. Related techniques were previously
described in \cite{kifmm04ying} and \cite{MR2005}.

Let $B_i$ denote a leaf box with associated index vector $I_{i}$.
We will describe the accelerated technique for constructing a matrix
$L_{i}$ and an index vector $I_{i}^{\rm rsk} \subset I_{i}$ such that
\begin{equation}
\label{eq:cup1}
A_{i}^{\rm row} = L_{i}\,A(I_{i}^{\rm rsk},:)
\end{equation}
holds to high precision. (The technique for finding $R_{i}$ and
$I_{i}^{\rm csk}$ such that (\ref{eq:fullID}) holds is analogous.)
For concreteness, suppose temporarily that the kernel $\mathcal{K}$
is the fundamental solution of the Laplace equation, $\mathcal{K}(x,y) = \frac{1}{2 \pi} \log|x-y|$.
The idea is to construct a small matrix $\tilde{A}_{i}^{\rm row}$
with the property that
\begin{equation}
\label{eq:cup2}
\mbox{Ran}\bigl(A_{i}^{\rm row}) \subseteq \mbox{Ran}\bigl(\tilde{A}_{i}^{\rm row}).
\end{equation}
In other words, the columns of $\tilde{A}_{i}^{\rm row}$ need to
span the columns of $A_{i}^{\rm row}$. Then compute an ID of the
\textit{small} matrix $\tilde{A}_{i}^{\rm row}$,
\begin{equation}
\label{eq:cup3}
\tilde{A}_{i}^{\rm row} = L_{i}\,\tilde{A}(I_{i}^{\rm rsk},:).
\end{equation}
Now (\ref{eq:cup2}) and (\ref{eq:cup3}) together imply that (\ref{eq:cup1}) holds.

\begin{figure}
\begin{center}
\begin{tabular}{cc}
\includegraphics[width=40mm]{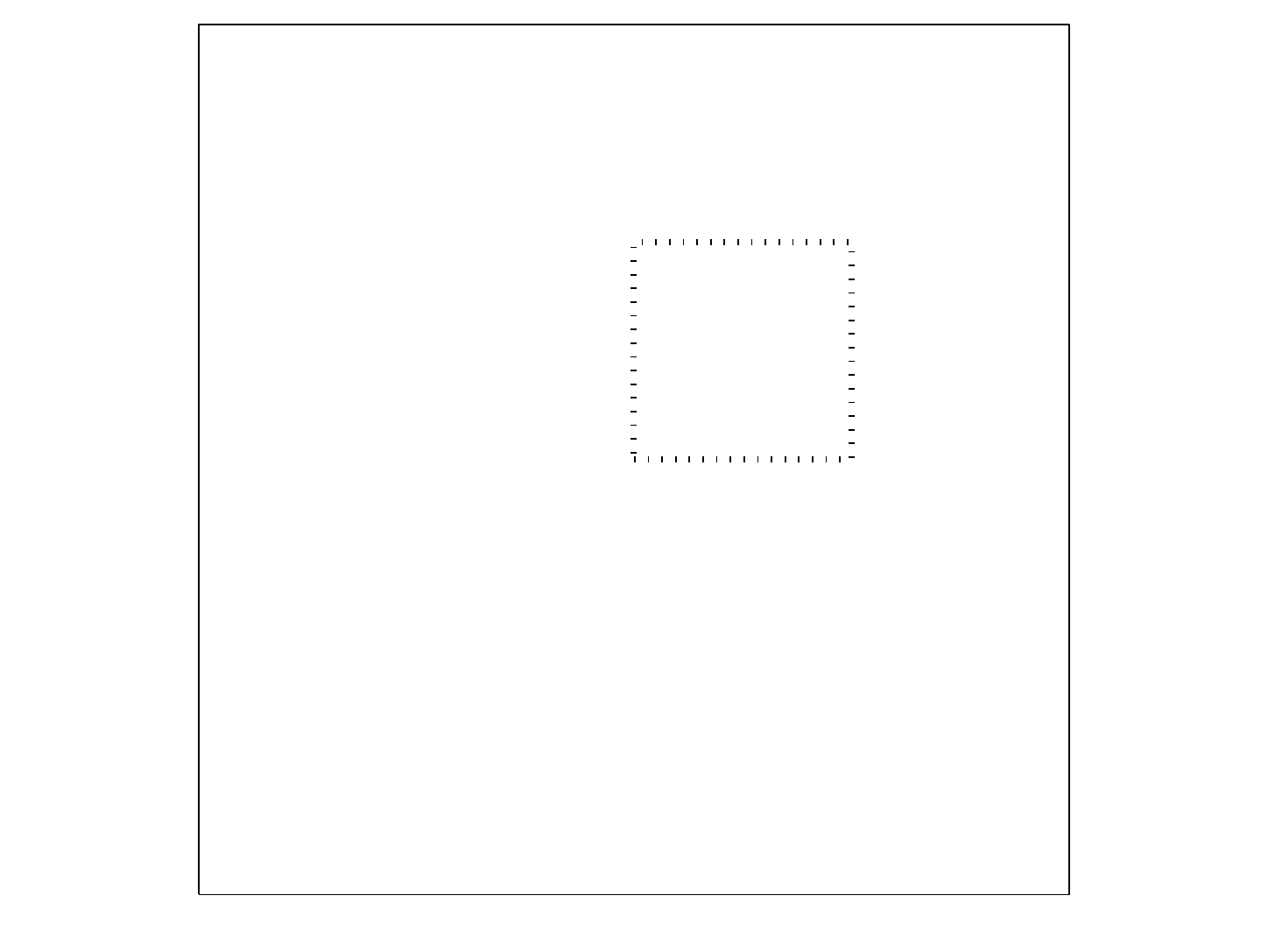}
&
\includegraphics[width=40mm]{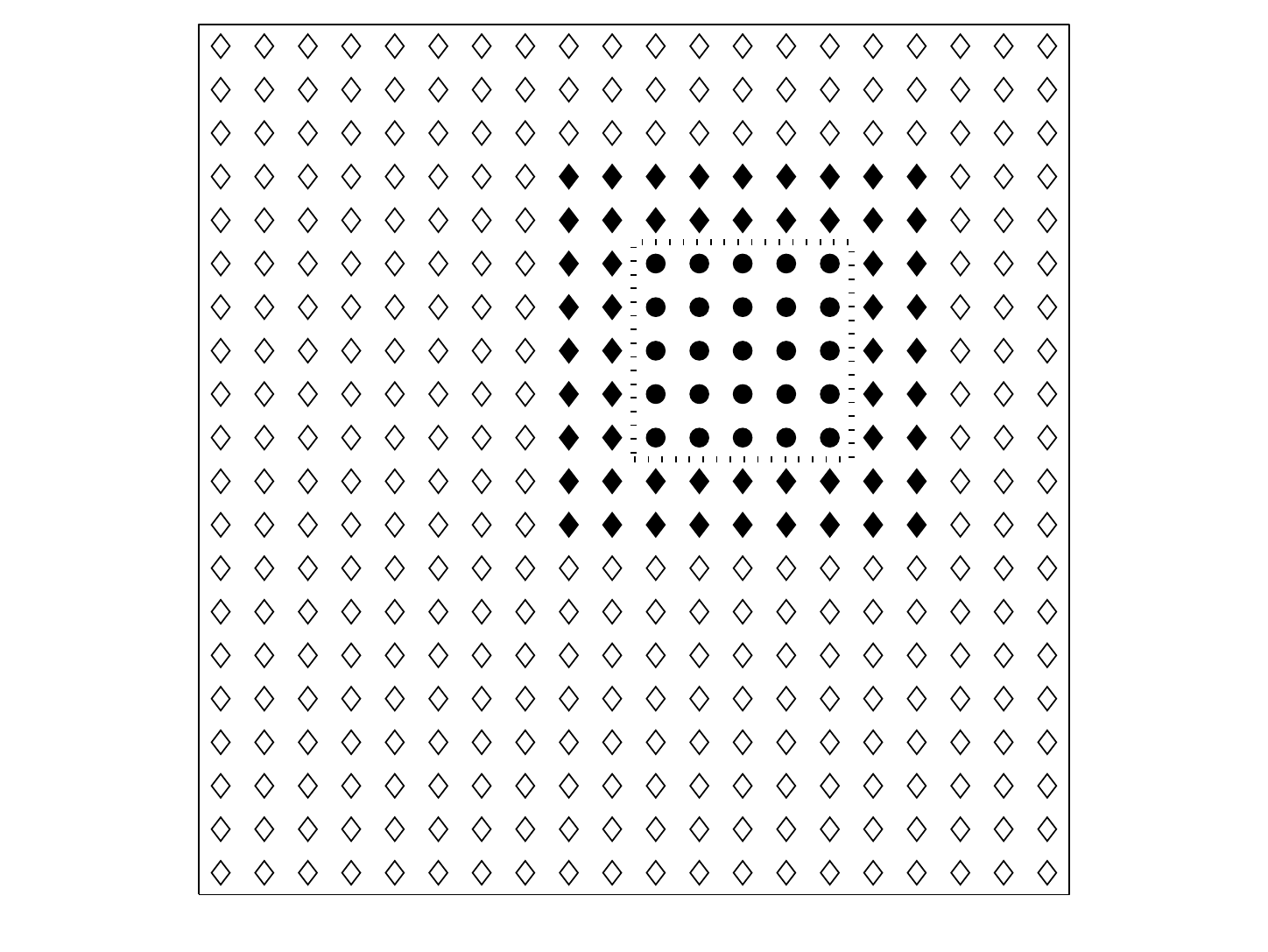} \\
(a) & (b)
\end{tabular}
\end{center}
\caption{(a) A domain $\Omega$ (solid) with a sub-domain $B_{i}$ (dotted).
(b) Target points in $B_{i}$ are circles, source points are diamonds,
and among the source points, the {\em proxy points} are solid.}
\label{fig:proxy}
\end{figure}

It remains to construct a small matrix $\tilde{A}_{i}^{\rm row}$
whose columns span the range of $A_{i}^{\rm row}$. To do this,
suppose that $v \in \mbox{Ran}(A_{i}^{\rm row})$, so that
for some vector $q \in \mathbb{R}^{N-m_{i}}$
$v = A_{i}^{\rm row}\,q$.
Physically, this means that the values of $v$ represent values
of a harmonic function generated by sources $q$ located outside
the box $B_{i}$. We now know from potential theory that any
harmonic function in $B_{i}$ can be replicated by a source density
on the boundary $\partial B_{i}$. The discrete analog of this statement
is that to very high precision, we can replicate the harmonic function
in $B_{i}$ by placing point charges in a thin layer of discretization
nodes surrounding $B_{i}$ (drawn as solid diamonds in Figure \ref{fig:proxy}(b)).
Let $\{z_{j}\}_{j=1}^{p_{i}}$ denote the locations of these points.
The claim is then that $v$ can be replicated by placing some
``equivalent charges'' at these points. In other words, we form
$\tilde{A}_{i}^{\rm row}$ as the matrix of size $m_{i} \times p_{i}$
whose entries take the form $\mathcal{K}(x_{r},z_{j})$ for $r \in I_{i}$,
and $j = 1,\,2,\,\dots,\,p_{i}$.

\begin{remark}
Figure \ref{fig:proxy} shows an example of how accelerated compression
works. Figure \ref{fig:proxy}(a) illustrates a domain $\Omega$ with 
a sub-domain $B_{i}$ (the dotted box). Suppose that $\varphi$ is a harmonic 
function on $B_{i}$. Then potential theory assures us that $\varphi$ can be 
generated by sources on $\partial B_{i}$, in other words $\varphi(x) = 
\int_{\Gamma}\mathcal{K}(x,y)\,\sigma(y)\,ds(y)$ for some density $\sigma$. 
The discrete analog of this statement is that to high precision, the harmonic 
function $\varphi$ can be generated by placing point charges on the proxy 
points $\{z_{j}\}_{j=1}^{p_{i}}$ marked with solid diamonds in Figure 
\ref{fig:proxy}(b). The practical consequence is that instead of
factoring the big matrix $A_{i}^{\rm row}$ which represents interactions between
all target points in $B_{i}$ (circles) and all source points (diamonds),
it is enough to factor the small matrix $\tilde{A}_{i}^{\rm row}$ representing
interactions between target points (circles) and proxy points (solid diamonds).
\end{remark}

\begin{remark}
The width of the layer of proxy points depends on the accuracy requested.
We found that for the Laplace kernel, a layer of width 1 leads to relative
accuracy about $10^{-5}$, and width 2 leads to relative accuracy $10^{-10}$.
For the Helmholtz kernel $\mathcal{K}(x,y) = H_{0}^{(1)}(\kappa|x-y|)$,
similar accuracy is \textit{typically} observed, but in this case, thicker
skeleton layers are recommended to avoid problems associated with resonances. 
(Recall that in classical potential theory, a solution to the Helmholtz equation may
require both monopole and dipole charges to be placed on $\partial B_{i}$.)
\end{remark}

\subsection{HSS matrix-vector multiplication}
\label{sec:HSS-multiply}
To describe the process of computing the inverse of an HSS matrix in
compressed form, it is convenient first to explain how matrix-vector
products can be computed.  The telescoping factorization \eqref{eq:telescope}  yields a fast matrix-vector apply
algorithm evaluating $u = A \sigma$. The structure of this algorithm
is similar to FMM (but simpler, as we do not treat the
near-field separately, approximating \emph{all} external interactions
with a single set of coefficients).  To emphasize the underlying
physical intuition, we refer to values $\sigma$ at
points as charges and to the values $u$ we want to compute as
potentials. On non-leaf boxes, we use notation $\phi^{\ell}_i$ for
the charges assigned to the skeleton points of the box, and
$u^{up,\ell}_i$ for computed potentials. The vector $\phi^{\ell}$ ($u^{\ell}$)
is the concatenations of all charges (potentials) of the boxes at level $\ell$.
At the finest level, we define $\phi^d := \sigma$.

\subsubsection*{Upward pass} The upward pass simply uses  the rectangular
block-diagonal matrices to compute $\phi^{\ell}$:   $\phi^{\ell} = R^{\ell+1} \phi^{\ell+1}$
Each block $R^{\ell+1}_i$ acts on the the subvector of charges
corresponding to the children of $B_i$ at level $\ell$.

\subsubsection*{Downward pass} We compute a potential $u_i$ for each box,
due to all \emph{outside} charges. starting from two boxes at level
$1$. For the top-level boxes $B_1$ and $B_2$, the values are computed 
directly, using the sibling interaction matrix forming $D^0$ as
defined by \eqref{eq:D-expressions}, in other words $u^1 = D^0 \phi^1$.

For boxes at the level $\ell > 1$, the  outside field on the boxes is
obtained as the sum of the fields interpolated to the boxes at level
$\ell$ using $L^\ell$ (``tall'' rectangular block-diagonal) and
contributions of the siblings through square block diagonal matrix
$D^\ell$:
\begin{equation*}
 u^{\ell} = D^{\ell} \phi^{\ell} + L^\ell u^{\ell-1}
\label{eq:dn}
\end{equation*}
At the leaf level, the last step is to transfer to the sample points
and add the field due to boxes themselves (self-interactions, stored
in the diagonal blocks $A_{i,i}$ of $D^{d}$), The actions of different
transformations are summarized in the following computational flow diagram:
\[
\xymatrix{
\sigma = \phi^{d} \ar[rr]^{R^d} \ar[rd]^{D^d} \ar[d]_{A^d \approx A}& & \phi^{d-1} \ar[rr]^{R^{d-1}} \ar[rd]^{D^{d-1}} \ar[d]_{A^{d-1}} & & \cdots   \ \ar[r]^{R^1} & \phi^{0} \ar[d]_{D^0 =A^0}\\
u = u^{d}  &\oplus \ar[l] & u^{d-1} \ar[l]^{L^{d}} &\oplus\ar[l]& \ar[l]^{L^{d-1}}\cdots&   u^{0} \ar[l]^{L^1}\\
}
\]

\subsection{Computing the HSS form of  $A^{-1}$}
\label{sec:HSS-inverse}
If $A$ is non-singular, the telescoping factorization \eqref{eq:telescope}
can typically be inverted directly, yielding an HSS representation of  $A^{-1}$,
which can be applied efficiently using the algorithm described in Section \ref{sec:HSS-multiply}.
The inversion process is also best understood in terms of the
variables introduced in Section~\ref{sec:HSS-multiply}.

First, we derive a formula for inversion of matrices having
single-level telescoping factorization $Z = F + LMR$.
(The matrix $Z$ on the first step coincides with $A = A^d$, and $F =
D^d$, but both $Z$ and $F$ are different  from $A^\ell$ and $D^\ell$
on subsequent steps).

We consider the system $Z \sigma = f$, and
define $\phi = \phi^{d-1} = R\sigma$,  $u = u^{d-1} =  M\phi$.
We then perform block-Gaussian elimination on the resulting block system:

\begin{equation}
\bbm F & L & 0 \\ -R & 0 & I \\ 0 & -I & M \ebm \bbm \sigma \\ u \\ \phi \ebm = \bbm f \\ 0 \\ 0 \ebm
\end{equation}

We form the auxiliary matrices $E = RF^{-1}L$, and $G = E + M$.  Then

\begin{equation}
\bbm F & 0 & 0 \\ 0 & E^{-1} & 0 \\ 0 & 0 & G \ebm \bbm \sigma \\ u \\ \phi \ebm = \bbm [I - LERF^{-1} + LEG^{-1}ERF^{-1}]f \\ [RF^{-1} - G^{-1}ERF^{-1}]f \\ ERF^{-1}f \ebm
\end{equation}
which yields the inverse of $Z$  by solving the block-diagonal system in the first line. Denoting $\tA[R] = ERF^{-1}$, $\tA[D] = F^{-1}(I - L\tA[R])$ and $\tA[L] = F^{-1}LE$, we obtain:

\begin{equation}
Z^{-1} = \tA[D] + \tA[L] (E+M)^{-1} \tA[R]
\label{eq:inv-recursion}
\end{equation}

We make a few observations:
\begin{itemize}
\item If $D,L$ and $R$ are block diagonal, then so are $E,\tA[R],\tA[L]$ and $\tA[D]$. 
This means that these matrices can be computed inexpensively via independent computations
that are local to each box.
\item The factors in the inverse can be interpreted as follows:
	\begin{itemize}
	\item $E^{-1} = RF^{-1}L$  can be viewed as a \emph{local solution operator} "reduced" to the set of skeleton points for each box. It maps fields to charge densities on these sets.
        \item $E + M$ maps charge densities on \emph{the union of skeleton points} to fields, adding diagonal $E\sigma$ and off-diagonal $M\sigma$ contributions.
	\end{itemize}
\end{itemize}

This inversion procedure can be applied recursively.
To obtain a recursion formula,  we define an auxiliary matrix
$\tA[A]^\ell = A^\ell + E^{\ell+1}$, for $\ell < d$, and $\tA[A]^d =
A^d = A$; $F^\ell = D^\ell + E^{\ell+1}$, and $F^d = D^d$.
Then  $\tA[A]^\ell$ satisfies  $ \tA[A]^\ell = F^\ell + L^\ell
A^{\ell-1} R^\ell$, and its inverse  by \eqref{eq:inv-recursion} satisfies

\begin{equation}
\left( \tA[A]^\ell\right)^{-1} = \tA[D]^\ell + \tA[L]^\ell (E^\ell +
A^{\ell-1} )^{-1} \tA[R]^\ell =
\tA[D]^\ell + \tA[L]^\ell (\tA[A]^{\ell-1} )^{-1} \tA[R]^\ell
\label{eq:inv-telescope}
\end{equation}
A fine-to-coarse procedure for computing the blocks of the inverse
immediately  follows from  \eqref{eq:inv-telescope}:
for each layer $\ell$, we first compute $F^\ell$, using $E^{\ell+1}$
for the finer layer (zero for the finest). $F^\ell$ determines
$\tA[D]^\ell$, $\tA[R]^\ell$ and $\tA[L]^\ell$, and $E^{\ell}$, to
be used at the next layer.

Algorithm~\ref{alg:HSS-inv} summarizes the HSS inversion algorithm;
it takes as input a tree $\cal T$ with index sets $I_i$
defined for leaf nodes, and matrices $D_i$, $R_i$ and $L_i$ for each
box $B_i$, and computes  components $\tA[D]_i$, $\tA[R]_i$
and $\tA[L]_i$ of the HSS form of  $A^{-1}$. This algorithm has
complexity $O(N^{3/2})$ for a volume integral equation in 2D.
It forms the starting point from which we derive an $O(N)$
algorithm in Section \ref{sec:inversion}.

\begin{algorithm}[h!]
\begin{center}
\begin{algorithmic}[1]
\FOR{ each   box $B_i$ in fine-to-coarse order}
\IF{$B_i$ is a leaf}
\STATE $F_i = D_i$
\ELSE
\STATE $F_i = D_i + \bbm  E_{c_1(i)} &\\ &E_{c_2(i)}  \ebm$
\ENDIF
\IF{$B_i$ is top-level}
\STATE $\tA[F]_i = F_i^{-1}$ \COMMENT{Direct inversion at the top level}
\ELSE
\STATE Compute $F_i^{-1}$
\STATE   $E_i = (R_i F_i^{-1} L_i )^{-1}$
\STATE $\tA[R]_i = E_i R_i F_i^{-1}$
\STATE $\tA[D]_i  = F_i^{-1} (I - L_i \tA[R]_i)$
\STATE $\tA[L]_i = F_i^{-1} L_i E_i$
\ENDIF
\ENDFOR
\end{algorithmic}
\end{center}
\caption{HSS matrix inversion}
\label{alg:HSS-inv}
\end{algorithm}

The hierarchical structure this algorithm computes is not entirely the same as that of 
matrix $A$: $\tA[L],\tA[R]$ are not interpolation matrices, and $\tA[D]$ has nonzero 
diagonal blocks. However, it can be converted to standard HSS format if needed via 
the simple re-formatting algorithm in \cite{GYMR2012, G2011thesis}.

\section{$O(N)$ Inverse Compression Algorithm}
\label{sec:inversion}

When applied to boundary integral equations (BIEs) in 2D, Algorithm 1
has optimal $O(N)$ asymptotic complexity for non-oscillatory 
kernels and a broad range of geometries. However, for volume 
integral equations in 2D, the typical complexity is $O(N^{3/2})$.

In this section, we first show why the algorithm
has higher complexity for volume problems, and then develop a faster 
algorithm in two steps. First, we modify the inverse algorithm to expose the
essential blocks of $\cA[A]^{-1}$, which can in turn be viewed as
compressible operators acting on 1D-like sets of points, and
avoid the need to use non-low-rank matrix-matrix products in the
algorithm. Then we demonstrate how these operators can be compressed
using one-dimensional HSS structures.

\subsection{Efficient skeleton construction in 2D}
\label{sec:skeleton_construction}

\subsubsection*{Skeleton size scaling}
\begin{figure}[t]
\begin{center}
\includegraphics[scale = 0.5]{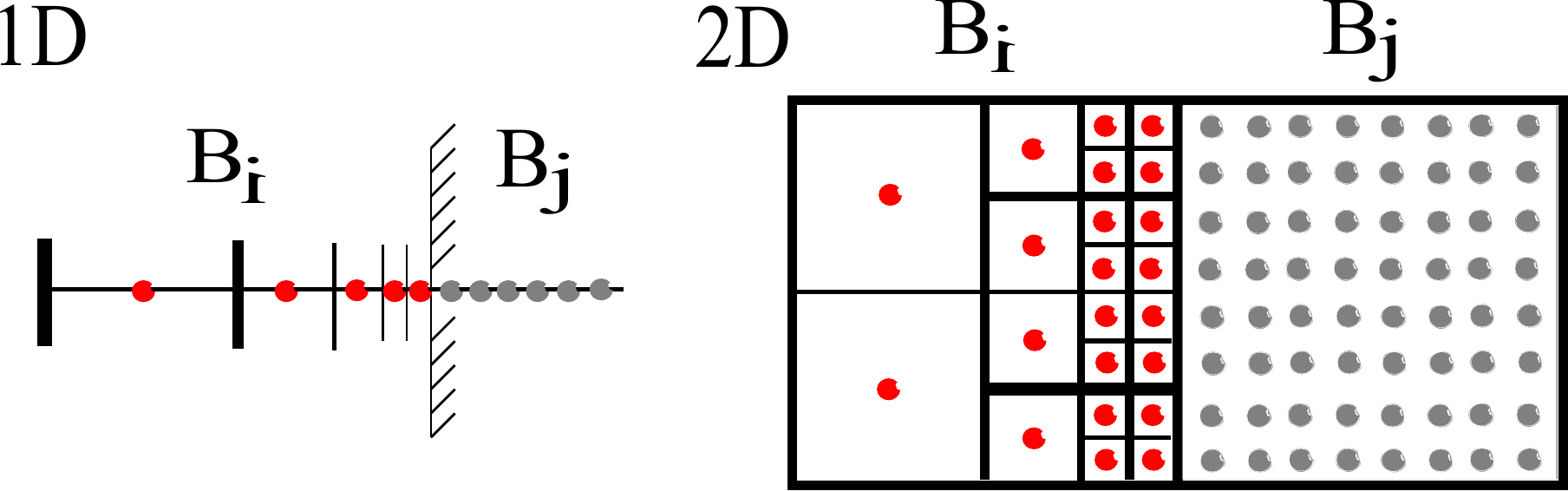}
\caption{ \textit{Interaction ranks in 1D and 2D. Source box $B_i$ is recursively subdivided into well-separated sets, whose interaction with $B_j$ is constant rank. This provides an upper bound for the overall interaction rank.} }
\end{center}
\label{fig:interaction_ranks}
\end{figure}

Let us first consider for a non-oscillatory kernel (e.g. Laplace) the rank of the interaction between two
neighboring boxes for dimensions $D = 1,2$. This example captures the essential behavior
and size of skeleton sets in different dimensions. (Figure \ref{fig:interaction_ranks})

Let $B_i$ be a box at level $\ell$ of the tree in $D$ dimensions, with
$m_i \sim n_{\ell} = N/2^{\ell}$ points.
We estimate the $\acc$-rank $k_i$ of
the interaction of $B_i$ with a box of the same size adjacent to
it. If sources and targets are well-separated
(the distance between the set and target points is at least the set's
diameter), the rank of their interaction matrix can be bounded by a
constant $p$ for a given $\epsilon$.  We perform a recursive
subdivision of our box $B_i$ into well-separated sets  until they
contain $p$ points of less. Then

\beqn
k_{i} \sim p \sum_{s=0}^{\log_2(n_{\ell}/p)/D} {2^{(D-1)s}}
\eeqn

For $D=1$, we have $\log_2(n_{\ell}/p)$ intervals, and so $k_i \sim
O(\log_2(n_{\ell}))$. For $D=2$,  we subdivide $\log_{4}(n_{\ell}/p)$ times in the
direction normal to the shared edge of the boxes, getting $2^{s}$
well-separated boxes. Then, $k_i \sim 2^{\log_{2}(n_{\ell}/p)/2} = O(n_{\ell}^{1/2})$
A simple calculation shows that this yields an estimate
$O(N)$ for the complexity of the 1D inversion algorithm,
(the logarithmic rank growth with box size does not affect the complexity,
because the number of boxes per level shrinks exponentially)
and  $O(N^{3/2})$ in the two-dimensional case.

\subsubsection*{Structure of skeleton sets in 2D}
The HSS hierarchical compression procedure
requires us to construct skeleton sets for
boxes at all levels.

The algorithm of Section~\ref{sec:hss} starts with constructing
skeletons for leaf boxes, using interpolative decomposition
of block rows with equivalent density acceleration.
For non-leaf boxes at level $\ell$, index sets are obtained by merging
skeleton sets of children and performing another
interpolative decomposition on the corresponding block row
of $A^\ell$ to obtain the skeleton.

This approach works for one-dimensional problems, but for two-dimensional ones applying
interpolative decomposition at all levels of the hierarchy is
prohibitively expensive even with equivalent density acceleration:
the complexity of each decomposition for a box $B_i$
is proportional to $k_i^3 \sim n_{\ell}^{3/2}$ in the two-dimensional
case. As a consequence, the overall complexity cannot be lower
than $O(N^{3/2})$.

\begin{figure}[t]
\begin{center}
\includegraphics[scale = 0.09]{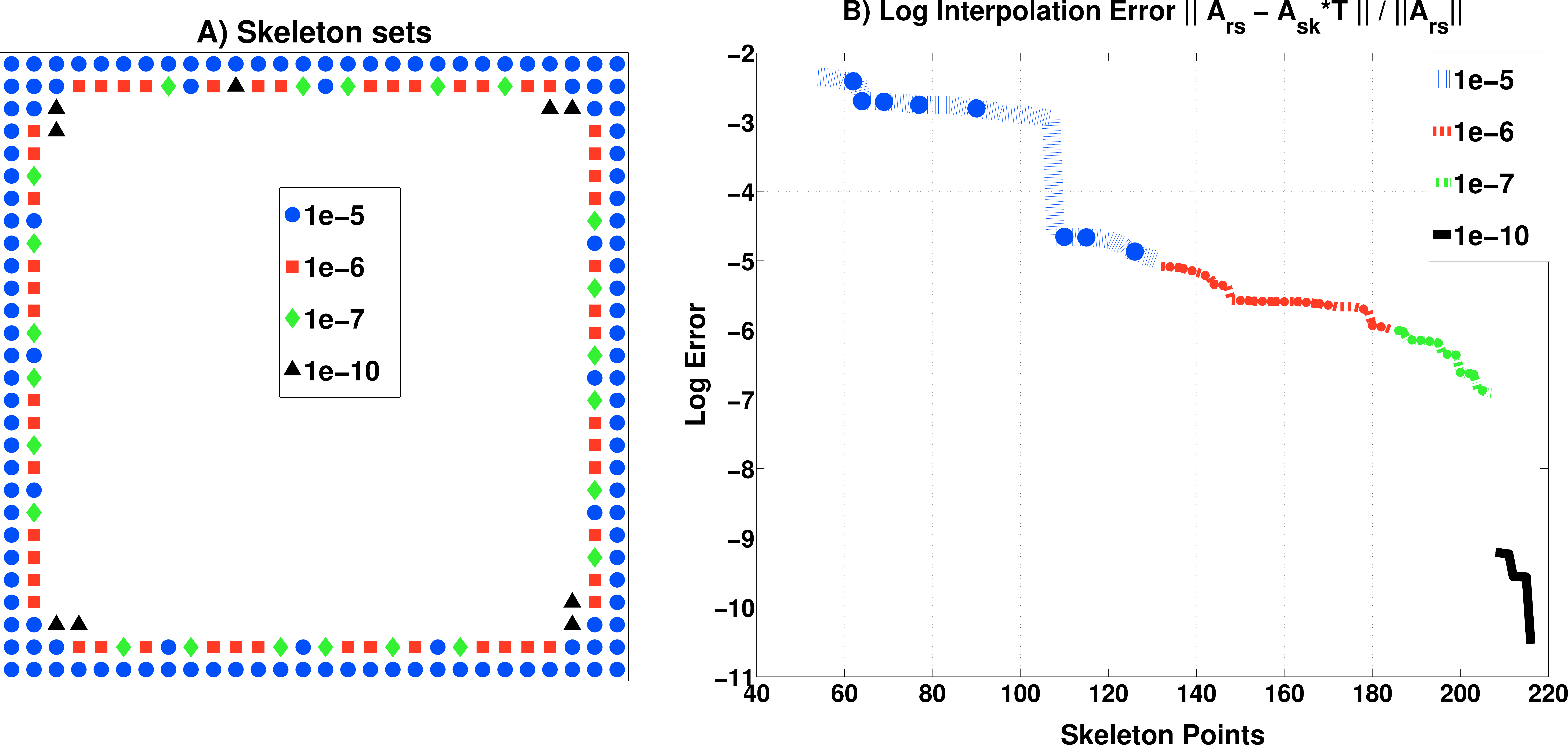}
\caption{ \textit{A) Skeleton sets picked by the interpolative decomposition for different relative interpolation accuracies B) Log plot of relative interpolation error depending on skeleton set size } }
\label{fig:laplace-skeletons}
\end{center}
\end{figure}

Our algorithm for constructing the skeleton sets at all levels
is based on the following crucial observation:
\emph{It is always possible to find an accurate set of skeleton points 
for a box by searching exclusively within a thin layer of points along
the boundary of the box.} 
This observation can be justified using
representation results from potential theory, cf.~Section \ref{sec:hss}.
It has also been substantiated by extensive numerical experiments.
Increasing target accuracy adds more points in deeper layers,
but the depth never grows too large:
for the kernels we have considered, the factorization selects one 
boundary layer for $\acc \sim 10^{-5}$, and two layers for $\acc \sim 10^{-10}$.
(Figure~\ref{fig:laplace-skeletons}).

This observation allows us to make two modifications to
the skeleton selection algorithm. First, we restrict the set of
points from which the skeletons are selected \emph{a priori}
to $m$ boundary layers. Second, rather than selecting
skeleton points for a parent box from the union of skeletons
of child boxes using an expensive interpolative decomposition,
we simply take all points in the boundary layers of the parent box.

More specifically,  $I_i = I_{c_{1}(i)}^{sk} \sqcup I_{c_{2}(i)}^{sk}$
is split into $I_i^{sk}$,  the skeleton of $B_i$,
consisting of all points of $I_i$ within $m$ layers of the boundary
of $B_i$, and $I_i^{rs} = I_i \setminus I_i^{sk}$ (residual index set),
consisting of points at the interface of two child boxes.
To obtain the interpolation operator
$T_i: I_i^{rs} \rightarrow I_i^{sk}$ provided by the
interpolative decomposition in the slower approach,
we use a proxy set $Z^{\rm proxy} = \{z_{j}\}_{j=1}^{p_{i}}$
(as described in Section \ref{sec:hss}), and compute $T_i$ from the following equation:
\beqn 
K(Z^{proxy},X_i^{sk})T = K(Z^{proxy},X_i^{rs}).
\eeqn

\begin{remark}.
The set from which the skeleton set is picked has fixed width.
This means that matrices acting on the skeleton set are compressible
(in the HSS sense) in a manner analogous to boundary integral operators
and admit linear complexity matrix algebra. These matrices in essence 
act like boundary-to-boundary operators on the box.
\end{remark}

\subsection{Overview of the modified inversion algorithm}
\label{sec:overview-inv}

Construction of a compressed inverse of an HSS matrix described in
Section~\ref{sec:HSS-inverse} proceeds in two stages:
first, the HSS form of the matrix $A$ is constructed, followed
by HSS inversion, which is performed without any additional
compression.  In our modified algorithm, significant changes
are made to the second stage, including compression of
all blocks. Not all blocks constructed at the first stage
(compression of $A$) are needed for the inverse construction,
so we reduce the first-stage algorithm to building the
tree and interpolation operators  $L_i,R_i$ only.

Examining  Algorithm~\ref{alg:HSS-inv},
we observe that all blocks formed for each box $B_i$ involve the
factors $E_i,F_i^{-1}$ and the interpolation matrices $L_i,R_i$.
One first step to save computation is to compute $\{E_i,F_i^{-1}\}$
only. Matrix-vector multiplication for blocks $\tA[D]$, $\tA[L]$ and
$\tA[R]$ needed for the inverse matvec algorithm can be implemented as
a sequence of matvecs for $L_i$, $R_i$, $E_i$ and $F_i^{-1}$.

The main operations in the algorithm are the two dense block inversions in lines $10-11$. Since $F_i$ is of size $m_i = k_{c_{1}(i)}+k_{c_{2}(i)}$ (merge of two skeleton sets) and $E_i$ is of size $k_i$, storage space of these blocks is $O(k_i^2) = O(n_{\ell})$ and inverting them costs $O(k_i^3)  = O(n_{\ell}^{3/2})$ floating point operations. This observation shows that it is  impossible to obtain linear complexity for HSS inversion if we store and invert $E_i$ and $F_i$ blocks densely, or even to build the interpolation matrices that form $L_i$ and $R_i$.

We first present a reformulation of the HSS inversion algorithm that partitions computation of $E_i$ and $F_i^{-1}$ into block operations, allowing us to compress blocks as low rank or one-dimensional HSS forms. As it is typically the case,  it is essential to avoid explicit construction of the blocks that we want compressed, i.e. they are constructed in a compressed form from the start.

\subsubsection*{Building and inverting $F_i$}
For a non-leaf box $B_i$, $F_i$ is a linear operator defined on $I_i = I_{c_{1}(i)}^{sk} \sqcup I_{c_{2}(i)}^{sk}$, that is, on the merge of skeleton points from its children. Using physical interpretation of Section~\ref{sec:HSS-inverse},
it maps charge distributions to fields on this set, adding contributions from local operators $E_{c_{1}(i)}, E_{c_{2}(i)}$ and sibling interactions.
We expect $F_i$ to have a rank structure similar to that of
$K[I_i,I_i]$.

Aside from needing a
compressed form and an efficient matvec for $F_i^{-1}$ to be used
in the inverse HSS matvec algorithm, we also use $F_i^{-1}$ to
construct $E_i^{-1} = R_i F_i^{-1} L_i$.
Let $\Pi_i$ be the permutation matrix that places skeleton points first.
 Matrices
$R_i = \bbm I \ T_i^{up} \ebm \Pi_i^T$ and $L_i = \Pi_i \bbm I \\
(T_i^{dn})^{T} \ebm$ have a block form, with sub-blocks
$T_i^{up}$ and $T_i^{dn}$, which we will show to be low-rank
(Section~\ref{sec:lowrank-interpol}). To construct $E_i$
efficiently, we need an explicit partition of $F_i$ into
blocks matching blocks of $R_i$ and $L_i$, i.e. corresponding
to skeleton index set $I^{sk}_i$ and residual index set $I^{rs}_i$ of $B_i$.
We use the following notation for the blocks of $F_i$:

\beqn  \Pi_i F_i  \Pi_i^{T} =
\bbm F_i[I_i^{sk},I_i^{sk}] & F_i[I_i^{sk},I_i^{rs}]  \\ F_i[I_i^{rs},I_i^{sk}]  & F_i[I_i^{rs},I_i^{rs}]  \ebm
= \bbm F_i^{sk} & F_i^{s \leftarrow r} \\ F_i^{r \leftarrow s} & F_i^{rs}, \ebm \eeqn
where we use $s$ to refer to the skeleton set of the parent
and $r$ to the ``residual'' set (the part of the union of the skeleton
set of the children not retained in the parent). 

To represent $F^{-1}$ we use these subblocks and perform block-Gaussian elimination. If $\Phi = \Pi_i F_i^{-1} \Pi^{T}$. Then:

\beqn \Phi_i =  \bbm  \phi_i^{sk} \ \phi_i^{s \leftarrow r} \\ \phi_i^{r \leftarrow s} \ \phi_i^{rs} \ebm \eeqn

We note that, by eliminating residual points first, $\phi_i^{sk}$ is the inverse of the Schur complement matrix $(S_i^{rs})^{-1} = [F_i^{sk} - F_i^{s\leftarrow r}(F_i^{rs})^{-1}F_i^{r \leftarrow s}]^{-1}$, and that only $\{ (F_i^{rs})^{-1},F_i^{r \leftarrow s},F_i^{s \leftarrow r},(S_i^{rs})^{-1} \}$ are needed to compute the blocks for the inverse $\Phi_i$. The routine in our inverse compression algorithm that builds $F_i^{-1}$ (Section~\ref{sec:buildFinv}) constructs these four blocks in compressed form.

\subsubsection*{Building $E_i$} Following the definition of $E_i$ and the
block structure of $F_i^{-1}$ we obtain the expression:
\beqn
E_i^{-1}
= \bbm I \ T_i^{up} \ebm
\bbm \phi_i^{sk} \ \phi_i^{s \leftarrow r}\\ \phi_i^{r \leftarrow s} \ \phi_i^{rs} \ebm
\bbm I \\ (T_i^{dn})^{T} \ebm
= \phi_i^{sk} + T_i^{up}\phi_i^{s \leftarrow r} + \phi_i^{r \leftarrow s} (T_i^{dn})^{T} + T_i^{up}\phi_i^{rs}(T_i^{dn})^{T}
\eeqn

As  $T_i^{up}$,$T_i^{dn}$ are low rank, the last three terms in this sum are low rank, too. Hence, $E_i$ can be computed as a low rank update of  $\phi_i^{sk}$, the inverse of $S_i^{rs}$ \\

We summarize the modified algorithm below; at this point it is a purely algebraic transformation, we explain in greater detail how the
new structure can be used to compress various matrices.
Next to each block we write between brackets the type of compression used in the $O(N)$ algorithm: \textbf{[LR]} (low rank) or \textbf{[HSS1D]}. The output of the new form consists of:

\begin{enumerate}
	\item The blocks of $F_i^{-1}$: $\{
	(F_i^{rs})^{-1},F_i^{r \leftarrow s},F_i^{s \leftarrow r},(S_i^{rs})^{-1} \}$
	\item The matrix $E_i$
\end{enumerate}
The entries of blocks of $F_i$ are evaluated using the formulas in Algorithm~\ref{alg:HSS-inv}.

\begin{algorithm}[h!]
\begin{center}
\begin{algorithmic}[1]
\FOR{ each   box $B_i$ in fine-to-coarse order}
\IF{$B_i$ is a leaf}
\STATE $F_i^{-1} = D_i^{-1} = K[I_i,I_i]^{-1}$
\ELSE
\STATE Compute blocks $F_i^{sk}$, $F_i^{rs}$ from $(E_{c_1(i)},E_{c_2(i)})$  \hfill \hfill     
 \linebreak In compressing and inverting the matrix $F_i$ for non-leaf boxes, we store sub-blocks of $F_i$ and apply block inversion formulae as explained above:

\STATE  $(F_i^{rs})^{-1} = (F[I_i^{rs},I_i^{rs}])^{-1}$ \hfill \textbf{[HSS1D]}
\STATE Compute blocks $F_i^{r \leftarrow s},F_i^{s \leftarrow r}$ \hfill \textbf{[LR]}
\STATE
$ (S_i^{rs})^{-1} = (F_i^{sk} - F_i^{s\leftarrow
r}(F_i^{rs})^{-1}F_i^{r \leftarrow s})^{-1}$ \hfill \textbf{[HSS1D]}
\ENDIF

\STATE
$E_i = \left( (S_i^{rs})^{-1} +  T_i^{up}\phi_i^{s \leftarrow r}
+ \phi_i^{r \leftarrow s}(T_i^{dn})^{T} +
T_i^{up}\phi_i^{rs}(T_i^{dn})^{T}\right)^{-1}$
\hfill \textbf{[HSS1D]}
\ENDFOR
\end{algorithmic}
\end{center}
\caption{Modified HSS inversion algorithm}
\label{alg:HSS-inv-mod}
\end{algorithm}
In the remaining part of this section we elaborate the details of efficient construction of all blocks using HSS and low-rank
operations.

\subsection{Compressed two-dimensional HSS inversion}
\label{sec:2dhss}

We use two compressed formats for various linear operators in the
algorithm: low-rank and one-dimensional HSS (dense-block HSS)
described in Section~\ref{sec:background}.
\subsubsection*{Operator Notation} To distinguish between linear operators
compressed in different ways, we use different fonts:
\begin{itemize}
  \item $X$ (in normal font) refers to an abstract linear operator, with no representation specified.

  \item $\cA[X]$ refers to a dense-block HSS representation of $X$;

  \item $\ttA[X]$ refers to a low-rank representation of $X$
\end{itemize}

The font used for interpolation operators like $R(:,J)=\bbm I \ T \ebm$ refers to the representation of $T$. Operations such as matrix-vector multiplies should be understood accordingly: e.g., $\ttA[X] v$ is evaluated using a low-rank factorization of $X$; if the rank is $q$, and vector size is $n$, the complexity is $O (qn)$. Similarly, $\cA[X] v$ is an $O (n)$ application of a dense-block HSS matrix. 

We say that the algorithms operating on per-box matrices
are fast if all operations involved have
cost and storage proportional to the block size $O(n^{1/2})$ or that
times a logarithmic factor $O(n^{1/2}\log^q(n))$. Our algorithm includes three main fast subroutines.

\begin{enumerate}
\item \textbf{INTER\_LOWRANK}: Interpolation matrices $\ttA[T]_i$ are built in low rank form.
\item \textbf{BUILD\_Finv}: Given dense-block HSS matrices $ \{ \cA[E]_{c_1(i)},\cA[E]_{c_2(i)} \}$, it computes $\{(\cA[F]_i^{rs})^{-1},\ttA[F]_i^{s \leftarrow r},\ttA[F]_i^{r \leftarrow s},(\cA[S]_i^{rs})^{-1} \}$ in their respective compressed forms.
\item \textbf{BUILD\_E}: Given
$\{(\cA[F]_i^{rs})^{-1},\ttA[F]_i^{s \leftarrow
r},\ttA[F]_i^{r \leftarrow s},(\cA[S]_i^{rs})^{-1} \}$ and
$\ttA[T]_i^{up},\ttA[T]_i^{dn}$, it computes $\cA[E]_i$ as a dense-block HSS matrix.
\end{enumerate}

\subsubsection{Fast Arithmetic}
\label{sec:fast_arithmetic}
There is a number of operations with low rank and HSS matrices which we must be able to perform efficiently:
\begin{itemize}
	\item \textbf{Dense-block HSS1D compression, inversion and matvec:} These are the algorithms in \cite{GYMR2012}, which we have outlined in Section 2, and we denote the corresponding routines as \textbf{HSS1D\_Compress} and \textbf{HSS1D\_Invert}. 
\item \textbf{\emph{Fast addition and manipulation of HSS1D matrices}}:
	\begin{itemize}
        	\item {\textbf{HSS1D\_Sum:}} Given two matrices $\cA[A]$ and $\cA[B]$ in HSS form, return an HSS form for $\cA[C] = \cA[A]+\cA[B]$
		\item {\textbf{HSS1D\_Split:}} Given a matrix $\cA[A]$ defined on $I_1  \cup I_2$, it produces the diagonal blocks $\cA[A]_1$ and $\cA[A]_2$ in HSS form.
		\item {\textbf{HSS1D\_Merge:}} Given matrices $\cA[A]_1$ and $\cA[A]_2$, it concatenates them to produce the block diagonal HSS matrix $\cA[A]$, and sorts its leaves accordingly.

	\item \textbf{\emph{Additional HSS compression routines}}:
	\begin{itemize}
		\item {\textbf{LR\_to\_HSS1D:}} Convert a low-rank
	operator to dense-block HSS form.
	        \item {\textbf{HSS1D\_Recompress:}} Using the algorithm by Xia \cite{xia2012complexity}, we re-compress an HSS 1D form to obtain optimal ranks. It is crucial to do so after performing fast arithmetic (e.g.~a sum) of HSS matrices.
\end{itemize}

	\end{itemize}
\end{itemize}

\begin{remark}
\textbf{Matrix-matrix products.} As we have mentioned before, a
        	feature of our algorithm is that it avoids using the
        	linear yet expensive matrix-matrix  product algorithm for
        	structured matrices.

We arrange the computations so that all matrix-matrix products are between a dense-block HSS $\cA[H]$ of size $k \times k$ and a low rank matrix $\ttA[K] = UV^{T}$ of rank $q$. Since $\cA[H]U$ requires only $q$ $O(k)$ fast matvecs, these products may be computed in $O(qk)$ work.
\label{rem:matrix-matrix}
\end{remark}

\begin{remark}
For every call to \textbf{HSS1D\_Sum} in the algorithms below, it
should be assumed that it is followed by a recompression step (a call
to \textbf{HSS1D\_Recompress}). Due to the theorem 5.3 of \cite{xia2012complexity}, this implies all dense-block 
HSS matrices presented are compact, that is, their blocks have a near-optimal size.
\end{remark}

\subsubsection*{Randomized interpolative decomposition \textbf{RAND\_ID}} The randomized sampling techniques in \cite{RandID07,RandID08,halko2011} speeding up the interpolative decomposition are critical to obtain the right complexity for the compression of low-rank operators. For an $m \times m$ dense-block  HSS matrix $A$ of rank $q$, the matvec has complexity $O(m)$ time. Assuming that at level $\ell$ of the binary 
tree $\cA[T]$, $m=O(n_{\ell}^{1/2})$ and $q=O(\log(n_{\ell}))$, the 
complexity of the randomized IDs is:

\beqn O(mq^2+q^3) = {O(n_{\ell}^{1/2}\log^2(n_{\ell}))} \eeqn

\subsubsection{Interpolation Operators in low-rank form}
\label{sec:lowrank-interpol}
We recall that interpolation operators $T_i$ are built with the binary
tree $\cA[T]$ using an interpolative decomposition, and that they are
inputs to the dense-block HSS inverse algorithms. They encode interactions between the residual points indexed as $I^{rs}_i$ and the exterior $X^{ext}_i$ as a linear combination of interactions with the skeleton points indexed by $I_i^{sk}$. In other words, they are solutions of the linear equation:

\beqn
K[X^{ext}_i,X_i^{sk}]T_i = K[X^{ext}_i,X_i^{rs}]
\eeqn

The acceleration proposed in Section \ref{sec:hss}  implies that interaction with outside points may be represented with a proxy $Z^{\rm proxy}_i$ of charges right outside the box boundary. We employ as many layers as are kept for $I_i^{sk}$, and remove points until $|Z^{\rm proxy}_i| = |X^{sk}_i| = k_i$. This yields an equation for  $T_i$:

\beqn
K[Z^{\rm proxy}_i,X_i^{sk}]T_i = K[Z^{\rm proxy}_i,X_i^{rs}]
\eeqn
$K[Z^{\rm proxy}_i,X_i^{sk}]$ encodes interactions between two close boundary layer 
curves, their distance being equal to the grid spacing $h$. It is invertible, 
ill-conditioned, and most importantly it has dense-block HSS structure.

$K[Z^{\rm proxy}_i,X_i^{rs}]$ encodes interactions with the interfacial points. 
Except for the closest layers, the interface is well-separated from the proxy, and so it is easy to see (via a multipole argument, or numerically as in Figure \ref{fig:inter-LR}) that this matrix is of logarithmic low rank:

\begin{figure}[H]
\begin{center}
~\includegraphics[scale = 0.09]{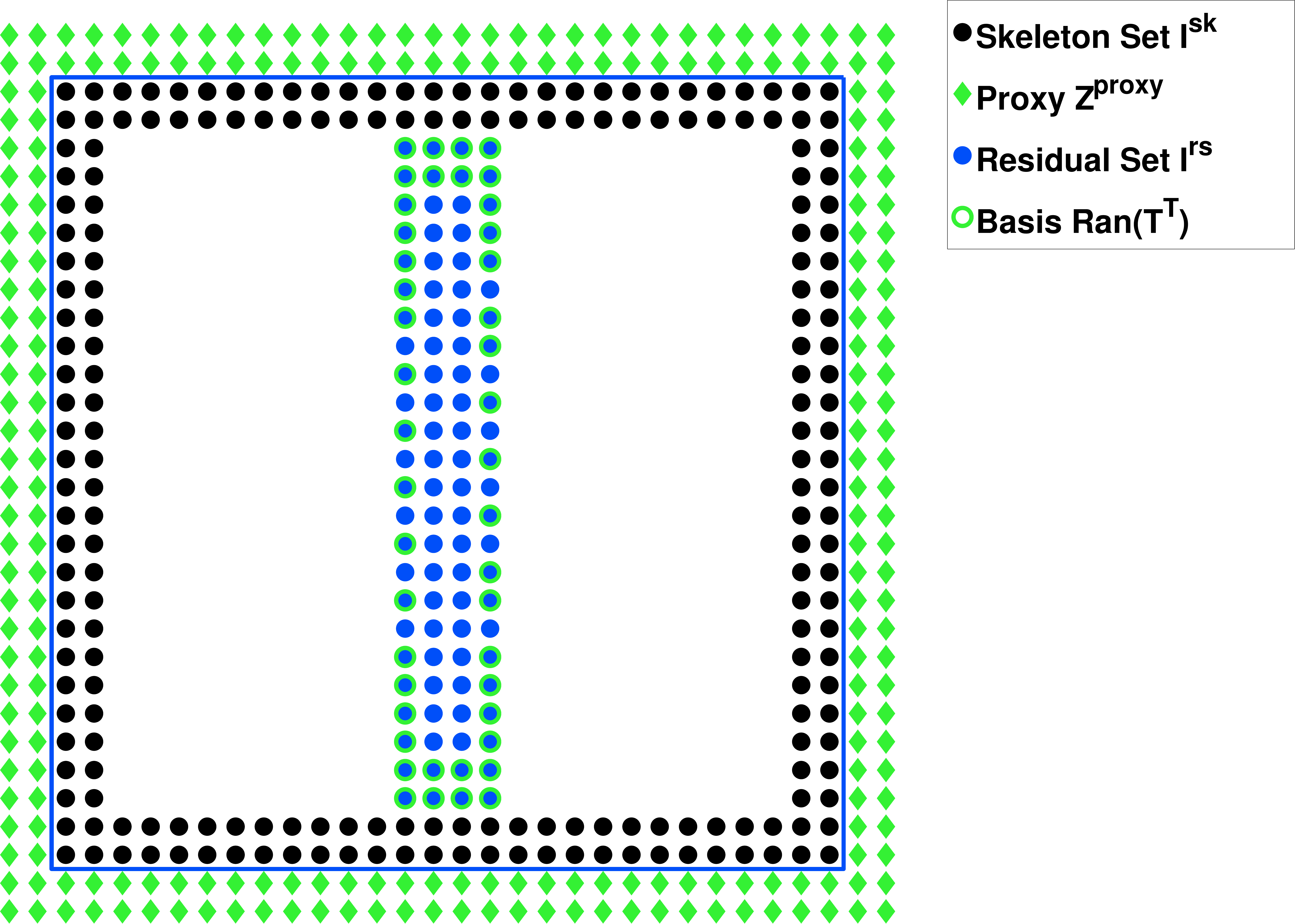}
\caption{\textit{Skeleton points are in black, residual points in blue and proxy points (diamonds) in green. We apply an ID with $\acc=1e-10$, and label subselected interface points also in green.}}
\label{fig:inter-LR}
\end{center}
\end{figure}

As a result, $T_i = K[Z^{\rm proxy}_i,X_i^{sk}]^{-1} K[Z^{\rm proxy}_i,X_i^{rs}]$ is also a low rank operator.
We describe a fast algorithm to compress this operator. We note that while the description above follows 
the case of boxes in the plane, these observations should hold in general, as long as most residual points 
are well separated from the proxy.

\begin{algorithm}[h]
\noindent\textbf{Input:} Box $B_i$ information, index sets $I^{sk}_i$ and $I^{rs}_i$;
\noindent\textbf{Output:} $\{U_{i},V_{i}\}$ such that $T_i = U_iV_i^{T}$.
\begin{center}
\begin{algorithmic}[1]

\STATE (i) \emph{Compress and Invert HSS 1D operator}
\STATE $\cA[K]^{s \rightarrow p} = $\textbf{HSS1D\_Compress}$(K,Z^{\rm proxy}_i,X_i^{sk},\acc)$
\STATE $(\cA[K]^{s \rightarrow p})^{-1} = $\textbf{HSS1D\_Invert}$(\cA[K]^{s \rightarrow p},\varepsilon)$
\STATE (ii) \emph{Randomized interpolatory decomposition}
\STATE $\cA[K]^{r \rightarrow p} = $\textbf{HSS1D\_Compress}$(K,Z^{\rm proxy}_i,X_i^{rs},\acc)$
\STATE $[T^{r \rightarrow p} , J^{r \rightarrow p}] =$\textbf{RAND\_ID}$(\cA[K]^{r \rightarrow p}, \acc);$
\STATE The ID gives us a low-rank decomposition of $K^{r \rightarrow p}$ (of $\acc$-rank q):
\STATE $U_i = (\cA[K]^{s \rightarrow p})^{-1}\ttA[K]^{r \rightarrow p}(:,J^{r \rightarrow p}(1:q));$
\STATE $V_i^{T}(:,J^{r \rightarrow p}) = \begin{bmatrix} I \ \ T^{r \rightarrow p} \end{bmatrix}$
\end{algorithmic}
\end{center}
\caption{Algorithm INTER\_LOWRANK}
\label{alg:inter-lowrank}
\end{algorithm}

\begin{remark}
In general, $Z_i^{proxy}$ will have slightly more points than $X_i^{sk}$, making $\cA[K]^{s \rightarrow p}$ a rectangular matrix. A fast HSS least squares algorithm such as in \cite{ho2012LS,dewildehierarchical} replaces the inversion and inverse apply in lines 3 and 8 of Algorithm \ref{alg:inter-lowrank}, with no impact in complexity. 
\end{remark}

\subsubsection{Compressed $F^{-1}$}
\label{sec:buildFinv}
We define some auxiliary index sets for the merge of two children boxes:
 $I_{i,c(i)}^{sk}$ are skeleton points shared with child $c(i)$ (on
 boundary layers) and  $I_{i,c(i)}^{rs}$ are residual children's skeleton points (in the middle interface). We define $I_i^{sk} = I_{i,c_1(i)}^{sk} \cup I_{i,c_2(i)}^{sk}$, and $I_i^{rs} = I_{i,c_1(i)}^{rs} \cup I_{i,c_2(i)}^{rs}$

\begin{figure}[H]
\begin{center}
\includegraphics[scale = 0.5]{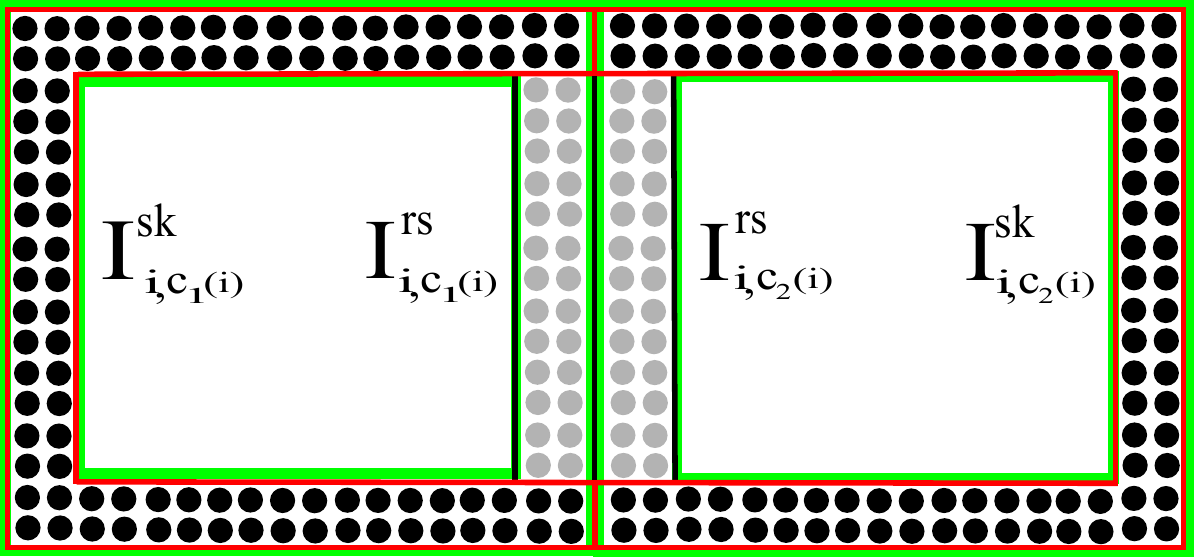}
\caption{\textit{Merge of two children's skeleton indices.}}
\end{center}
\end{figure}

We recall that first, matrix $F$ is built and separated into blocks
corresponding to skeleton $I_i^{sk}$ and residual $I_i^{rs}$
sets. These two sets can be seen as arranged  along one-dimensional
curves (boundary of the parent box and interface between siblings).
We order $I_i^{rs}$ in the direction of the interface, and $I_i^{sk}$ cyclically around the box boundary.

Then  $\cA[F]_i^{rs}=F(I_i^{rs},I_i^{rs})$ and  $\cA[F]_i^{sk} =
F(I_i^{sk},I_i^{sk})$  can be constructed in  dense-block HSS form.

The two off-diagonal blocks $\{ \ttA[F]_i^{s\leftarrow r},\ttA[F]_i^{r\leftarrow s} \}$ encode interaction between sets that are close only at few points, and can thus be compressed as low-rank operators.

\subsubsection*{\textbf{BUILD\_Finv} routine}
To clarify the structure of the construction of $F^{-1}$, in Algorithm 4 we show this routine in two ways. On the left, we indicate the original dense computation, and on which line of the reformulated algorithm in Section 2.5 it occurs. On the right, we indicate the set of fast operations that replaces it in the linear complexity algorithm.

\begin{algorithm}[h!]
\noindent\textbf{Input.}
Children's matrices $\cA[E]_{c_1(i)}$ and $\cA[E]_{c_2(i)}$ in HSS form (each defined on skeletons $I_{c_j(i)}^{sk}$). \\
\noindent\textbf{Output.}
\begin{enumerate}
\item $(\cA[F]_i^{rs})^{-1}$ in dense-block HSS form on the residual set $I^{rs}_i$.
\item ${\ttA[F]_i^{s \leftarrow r},\ttA[F]_i^{r \leftarrow s}}$ as low-rank operators.
\item Inverse of Schur complement matrix $(\cA[S]_i^{rs})^{-1}$ in
dense-block HSS form on $I^{sk}_i$.
\end{enumerate}
\noindent
\begin{tabular}{|l|l|l|}
\hline
\multicolumn{3}{|l|}{(i) Obtain diagonal blocks of E by splitting into boundary and interface}
\\
\hline
$E^{sk}_{c_j(i)} = E_{c_j(i)}[I^{sk}_{i,c_j(i)},I^{sk}_{i,c_j(i)}]$ && $[\cA[E]_{c_j(i)}^{sk},\cA[E]_{c_j(i)}^{rs}] = $ \textbf{HSS1D\_Split}$(\cA[E]_{c_j(i)},I_{i,c_j(i)}^{sk})$ \\
$E^{rs}_{c_j(i)} = E_{c_j(i)}[I^{rs}_{i,c_j(i)},I^{rs}_{i,c_j(i)}]$ && \\
\hline
\multicolumn{3}{|l|}{(ii) Build diagonal blocks of F: merge diagonal blocks of E, compress off-diagonal blocks}
\\
\hline
Line 5: && $\cA[E]_i^{sk,dg} =$ \textbf{HSS1D\_Merge}$(\cA[E]_{c_1(i)}^{sk},\cA[E]_{c_2(i)}^{sk});$  \\
$F^{sk} = \bbm E_{c_1}^{sk} & 0 \\ 0 & E_{c_2}^{sk} \ebm + \bbm 0 & K[I_{c_1}^{sk},I_{c_2}^{sk}] \\ K[I_{c_2}^{sk},I_{c_1}^{sk}] & 0 \ebm$ &&
      $\cA[K]_i^{sk,off} =$ \textbf{HSS1D\_Compress}$(K,I_{i}^{sk});$ \\
&& $\cA[F]_i^{sk}$ = \textbf{HSS1D\_Sum}$(\cA[E]_i^{sk,dg},\cA[K]_i^{sk,offd});$ \\
\hline
Line 5: && $\cA[E]_i^{rs,dg} =$ \textbf{HSS1D\_Merge}$(\cA[E]_{c_1(i)}^{rs},\cA[E]_{c_2(i)}^{rs});$  \\
$F^{rs} = \bbm E_{c_1}^{rs} & 0 \\ 0 & E_{c_2}^{rs} \ebm + \bbm 0 & K[I_{c_1}^{rs},I_{c_2}^{rs}] \\ K[I_{c_2}^{rs},I_{c_1}^{rs}] & 0 \ebm$ &&
      $\cA[K]_i^{rs,off} =$ \textbf{HSS1D\_Compress}$(K,I_{i}^{rs});$ \\
&& $\cA[F]_i^{rs}$ = \textbf{HSS1D\_Sum}$(\cA[E]_i^{rs,dg},\cA[K]_i^{rs,offd});$ \\
\hline
\multicolumn{3}{|l|}{(iii) Inverse of $F^{rs}_i$:}
\\
\hline
Line 6: && \\
$(F^{rs}_i)^{-1} = F[I^{rs}_i,I^{rs}_i]^{-1}$ && $(\cA[F]_i^{rs})^{-1}$ = \textbf{HSS1D\_Invert}$(\cA[F]_i^{rs});$ \\
\hline
\multicolumn{3}{|l|}{(iv) Low rank decompositions for $\ttA[F]_i^{s  \leftarrow r}$ and	$\ttA[F]_i^{r \leftarrow s}$ using Randomized IDs:}
\\
\hline
Line 7: && \\
$\ttA[F]_i^{s  \leftarrow r} = F[I^{sk}_i,I^{rs}_i]$ && $[T_i^{r},J_i^{r}] =$\textbf{RAND\_ID}$(F_i^{s \leftarrow r} , \acc);$ \\
$\ttA[F]_i^{r  \leftarrow s} = F[I^{rs}_i,I^{sk}_i]$ && $[T_i^{s},J_i^{s}] =$\textbf{RAND\_ID}$(F_i^{r \leftarrow s} , \acc);$ \\
\hline
\multicolumn{3}{|l|}{(v) Schur complement as a low rank perturbation of $\cA[F]_{i}^{sk}$:}
\\
\hline
Line 8: && $\cA[P]_i^{sk}$ = \textbf{LR\_to\_HSS1D}$( - \ttA[F]_i^{s \leftarrow r}(\cA[F]_i^{rs})^{-1}\ttA[F]_i^{r \leftarrow s});$ \\
$(S^{rs})^{-1} = [F[I_i^{sk},I_i^{sk}] - F^{s\leftarrow r}(F^{rs})^{-1}F^{r \leftarrow s}]^{-1}$ &&
$\cA[S]_i^{rs} =$ \textbf{HSS1D\_Sum}$(\cA[F]_{i}^{sk},\cA[P]_{i}^{sk});$ \\
&& $(\cA[S]_i^{rs})^{-1} =$ \textbf{HSS1D\_Invert}$(\cA[S]_i^{rs});$ \\
\hline
\end{tabular}
\caption{BUILD\_Finv}
\label{alg:build-inv}
\end{algorithm}

\subsubsection*{Applying $F^{-1}$ to a vector} Once we have the four
compressed blocks $\{ (\cA[F]_i^{rs})^{-1},  \ttA[F]_i^{s \leftarrow
r},\ttA[F]_i^{r \leftarrow s}, (\cA[S]_i^{rs})^{-1} \}$, the
routine \textbf{APPLY\_Finv} can be used to compute the fast product
$\sigma=F_i^{-1}u$.  This is done in a straightforward way  splitting
$u$ into $[ u[I_i^{sk}] , u[I_i^{rs}]]$ and using  fast matvecs of the
blocks of $F^{-1}$.

\subsubsection{Compressed $\cA[E]$}
\label{sec:buildE}
The last piece that is required is a routine that constructs $\cA[E]$ as a dense-block HSS matrix from $F^{-1}$ in  compressed form. For a box $B_i$, $E_i$ is a linear operator defined on the skeleton set $I_i^{sk}$. Its inverse is obtained by applying $F_i^{-1}$ to this set using interpolation operators:

\beqn E_i^{-1} = L_i F_i^{-1} R_i \eeqn

Hence, the physical intuition again is that $E_i$ has a rank structure similar to that of $K[I_i^{sk},I_i^{sk}]$. In \cite{bremer2011fast,MR07scat,chen2002} its inverse is called a \emph{reduced scattering matrix}. From a purely algebraic point of view, we observe that $\cA[E]$ is a low-rank perturbation of $\cA[S]^{rs}$. Hence, if $\cA[S]^{rs}$ has
dense-block HSS structure, so does $\cA[E]$.

\subsubsection*{BUILD\_E routine}
Recalling the definition of $E_i$ in  Algorithm~\ref{alg:HSS-inv-mod},
\beqn
E_i^{-1} = R_i F_i^{-1} L_i = \bbm I \ \ttA[T]_i^{up} \ebm
\bbm \phi_i^{sk} \ \phi_i^{s \leftarrow r}\\ \phi_i^{r \leftarrow s} \ \phi_i^{rs} \ebm
\bbm I \\ (\ttA[T]_i^{dn})^{T} \ebm
= \phi_i^{sk} + \ttA[T]_i^{up}\phi_i^{s \leftarrow r} + \phi_i^{r \leftarrow s}(\ttA[T]_i^{dn})^{T} + \ttA[T]_i^{up}\phi_i^{rs}(\ttA[T]_i^{dn})^{T}
\eeqn
we observe the last three matrices are low-rank. Using the Schur
complement formulae we obtain explicit factorizations for each one,
recompress as a single low-rank matrix and then convert it to the
dense-block HSS  form using \textbf{LR\_to\_HSS1D}.

\begin{algorithm}[h!]
\noindent\textbf{Input.}
$F^{-1}$ as blocks:
$(\cA[F]_i^{rs})^{-1},\ttA[F]_i^{s \leftarrow
r},\ttA[F]_i^{r \leftarrow s}$,$(\cA[S]_i^{rs})^{-1}$;
$L_i$ and $R_i$ in low-rank form:
$\ttA[T]_i^{up}$,$J_i^{up}$,$\ttA[T]_i^{dn}$, and $J_i^{dn}$,
$\ttA[T]_i^{up} = U_i^{up}V_i^{up}$, $(\ttA[T]_i^{dn})^{T} = U_i^{dn}V_i^{dn}$.
\noindent\textbf{Output.} $\cA[E]_i$ in HSS form (on the curve $I_i^{sk}$)

\begin{center}
\begin{algorithmic}[1]
\STATE $V_{i,1}^{E} = V_i^{up}\phi^{s \leftarrow r}; \ U_{i,1}^{E} = U_i^{up}$
\STATE $U_{i,2}^{E} = \phi^{r \leftarrow s}U_i^{dn}; \ V_{i,2}^{E} = V_i^{dn}$
\STATE $U_{i,3}^{E} = U_i^{up}V_i^{up}\phi^{rs}U_i^{dn}; \ V_{i,3}^{E} = V_i^{dn}$
\STATE Recompress: $U_i^E(V_i^E)^{T} = \sum_{p=1}^{3} {U_{i,p}^{E}(V_{i,p}^{E})^{T}}$
\STATE $\cA[M]_{i} = \textbf{LR\_to\_HSS1D}(U_{i}^{E},V_{i}^{E})$
\STATE $\cA[E]_i^{-1} = \textbf{HSS1D\_Sum}((\cA[S]_i^{rs})^{-1} , \cA[M]_{i})$
\STATE $\cA[E]_i = \textbf{HSS1D\_Invert}(\cA[E]_i^{-1})$
\end{algorithmic}
\end{center}
\caption{BUILD\_E}
\label{alg:builde}
\end{algorithm}

We can perform the low-rank matrix products (of rank $q_i$) in
Algorithm~\ref{alg:builde} in $O(m_i q_i)$ or $O(m_i q_i^2)$ operations, since all of the products involved
are fast ($O(k_i)$ or $O(m_i-k_i)$, where $|I_i| = m_i$ and $|I_i^{sk}| = k_i$).

\subsection{Inverse matrix-vector multiplication}
Once the compressed form of the inverse is obtained, it can be
efficiently  applied to a right-hand side vectors 
with the algorithm of  Section~\ref{sec:HSS-multiply},
but using fast algorithms for our compressed representations of blocks.
We present it here for completeness  (Algorithm~\ref{alg:inv-matvec}).

\begin{algorithm}[h!]
\noindent\textbf{Input}
\begin{enumerate}
  \item The binary tree $\mathcal{T}$ including skeleton set indices and $R_i,L_i$
  associated with boxes;
  \item HSS-compressed inverse: per-box blocks forming $F_i$ and $\cA[E]_i$.
  \item the vector $f$ of field values defined at source points.
\end{enumerate}
\subsubsection*{Output} The output is $\sigma = A^{-1}f$, where $A^{-1}$ is the compressed inverse.

\begin{center}
\begin{algorithmic}[1]
\COMMENT {Upward Pass, compute $u_i^{up}$}
\FOR{ each box $B_i$ in fine-to-coarse order}
\IF{$B_i$ is a leaf}
\STATE $u = f(I_i^{up})$
\ELSE
\STATE $u = \bbm u_{c_1(i)}^{up} \\ u_{c_2(i)}^{up} \ebm $
\ENDIF
\STATE Multiplication by $\tA[R]_i = E_i R_i F_i^{-1}$:
\STATE $\varphi_i = \textbf{APPLY\_Finv}(u)$
\STATE $u_i^{up} =  \cA[E]_i \ttA[R]_i \varphi_i$

\ENDFOR
\STATE \COMMENT {Downward Pass: compute $\phi^{dn}_i$; the result $A^{-1}f$ on leaf $B_i$ is  $\phi^{dn}_i$}
\STATE $\phi^{dn}_{top} = 0$ \ , \ $\nu_{top}^{dn} = u$
\FOR{ each $B_i$ in coarse-to-fine order}

\STATE Multiplication by $\tA[D]_i = F_i^{-1}[I - L_i\tA[R]_i]$
\STATE Define $u$ as above
\STATE $\nu_i^{dn} = u - L_i u_i^{up} $

\STATE Multiplication by $\tA[L]_i = F_i^{-1} L_i E_i$:
$\rho_i^{dn} = \ttA[L]_i \cA[E]_i \phi_i^{dn}$

\STATE Add both contributions multiplying by the common factor $F_i^{-1}$:
\IF {$B_i$ is a leaf}
\STATE $ \sigma(I_i^{dn}) = \textbf{APPLY\_Finv}(\nu_i^{dn} + \rho_i^{dn})$
\ELSE
\STATE $ \bbm \phi^{dn}_{c_1(i)} \\ \phi^{dn}_{c_2(i)} \ebm = \textbf{APPLY\_Finv}(\nu_i^{dn} + \rho_i^{dn}) $
\ENDIF
\ENDFOR
\end{algorithmic}
\end{center}
\caption{HSS inverse matrix-vector multiplication.}
\label{alg:inv-matvec}
\end{algorithm}

We notice that application of $\tA[R]_i, \tA[D]_i$ and $\tA[L]_i$ is
substituted by fast matrix vector multiplication of
$\ttA[L]_i, \ttA[R]_i$ , $\cA[E]_i$ and the blocks comprising
$F_i^{-1}$. All of these have complexity $O(n_{\ell}^{1/2})$ or $O(n_{\ell}^{1/2}\log(n_{\ell}))$.

\section{Complexity Estimates}
\label{sec:complexity}
In this section we estimate the computational complexity of the
dense-block and the compressed-block algorithms presented in
Section \ref{sec:inversion}, under a number of assumptions on the
rank structure of blocks of $\cA[A]$.

We first define a framework to estimate work and storage for
algorithms defined on a binary tree $\cA[T]$, and use it to analyze
all types of algorithms described in this paper.
We consider both non-translation invariant  (\textbf{NTI}), and
translation-invariant kernels (\textbf{TI})  for which significant
performance gains can be obtained.

\subsection{Complexity of algorithms on binary tries}
\label{sec:framework}

All our algorithms compute and store matrix blocks associated with boxes
organized into a binary tree $\cA[T]$. To produce complexity estimates
for a given accuracy $\acc$, we introduce bounds for work $W_{\ell}(n_{\ell},\acc)$ and storage $M_{\ell}(n_{\ell},\acc) $ at each level $\ell$ of the tree, where $n_{\ell}=2^{-\ell}N$ is the maximum number of points in a box at this level (we assume that the work per box on a given level has small variance).

\begin{lemma}
\label{lemma:complexity}
 Let $n_{\ell}=2^{-\ell}N$, and $d = \log_2(N/n_{max})$, and exponents
 $p,q \ge 0$. Then if  $W_{\ell}(n_{\ell},\acc)$ has the form $C_{\acc} n_{\ell}^p \log_2^q(n_{\ell})$,
the total work has complexity

\beqn
\textbf{NTI: }\sum_{\ell = 0}^{d} {2^{\ell} W_{\ell}(n_{\ell},\acc)} =
\left\{ \begin{array}{lr}
O(N) & : 0 \le p<1 \\
O(N \log_2^{q+1}N) & : p = 1 \\
O(N^p \log_2^{q}N) & : p > 1
\end{array} \right.
\eeqn

\beqn
\textbf{TI: } \sum_{\ell = 0}^{d} {W_{\ell}(n_{\ell},\acc)} = O(N^p \log^{q}N)
\eeqn
\end{lemma}

For NTI algorithms, on each level $\ell$ we obtain an estimate by adding the bound for work  per box for the $2^{\ell}$ boxes. The polynomial growth in $W_{\ell}$ or $M_{\ell}$ is compensated by the fact that the number of boxes decreases exponentially going up the tree.

If  the rate of growth of $W_{\ell}$ is slower than linear ($p<1$),
 the overall complexity is linear. If  $W_{\ell}$  grows linearly,  we
 accumulate a $\log N$ factor going up the tree.
If the growth of  $W_{\ell}$ is  superlinear ($p>1$), the work
 performed on the top boxes dominates, and we obtain the same
 complexity as in $W_{\ell}$ for the overall algorithm.

In the TI case, work/storage on the top boxes dominates the
calculation, since only one set of matrices needs to be computed and
stored per level.  Hence, the interpretation is simpler: the
single-box bound  for complexity at the top levels reflects the overall complexity.

\subsection{Assumptions on matrix structure}
\label{sec:assumptions}

Complexity estimates for the work per box require assumptions on  
matrix structure. Let us restate some assumptions already
made in Sections \ref{sec:background} and \ref{sec:inversion}:

\begin{enumerate}

  \item \label{asm:skeleton} \textbf{Skeleton size scaling:} 
  The maximum size of skeleton sets for boxes at level $\ell$ 
  grows as $O(n_{\ell}^{1/2})$. This determines the size of 
  blocks within the HSS structure, and can be proved for
  non-oscillatory PDE kernels in 2D.

  \item \label{asm:eqdensity} \textbf{Localization:} Equivalent 
  densities may be used to represent long range interactions to
  within any specified accuracy, cf.~Section \ref{sec:hss}.
    
  \item \label{asm:skeleton-struct} \textbf{Skeleton structure:} 
  The skeleton set for any box may be chosen from within a thin 
  layer of points close to the boundary of the box, cf.~Section \ref{sec:skeleton_construction}.  
  
  \item \label{asm:compressed-blocks} \textbf{Compressed block structure:} 
  Experimental evidence and physical intuition from scattering problems
  allows us to assume that the blocks of $F$ and $E$ discussed in 
  Section~\ref{sec:inversion} have one-dimensional HSS or low-rank 
  structure, with logarithmic rank growth ( $O(\log(n_{\ell}))$ ).
\end{enumerate}

\noindent
These assumptions arise naturally in the context of
solving integral equations with non-oscillatory PDE kernels in 2D. All
assumptions excluding the last one are relevant for both dense and
compressed block algorithms. The last one is needed only for the
compressed-block algorithms.

We note assumption \ref{asm:skeleton-struct} implies \ref{asm:skeleton}: being able to pick skeletons from a thin boundary layer determines how their sizes scale. We mention them separately to distinguish their roles on the design and complexity analysis of our algorithms: while \ref{asm:skeleton} mainly impacts block sizes on the outer HSS structure, \ref{asm:skeleton-struct} is much more specific and refers to a priori knowledge of skeleton set structure which we exploit extensively in the compressed-block algorithm. 

\subsection{Estimates}
\label{sec:estimates}

We analyze work and storage for the algorithms of  Section~\ref{sec:inversion}. Since they all use the same 
set of fast operations (see Section \ref{sec:fast_arithmetic}), we can make unifying observations:

\subsubsection*{Work}
\begin{enumerate}
\item Assumptions \ref{asm:skeleton}  and \ref{asm:skeleton-struct}
imply that our fast subroutines perform operations with HSS blocks of
size $O(n_{\ell}^{1/2})$.
Further, Assumption~\ref{asm:compressed-blocks} states that these
behave like HSS matrices representing boundary integral operators, 
for which all one-dimensional HSS operations are known to be linear 
in matrix size. Thus, all \textbf{HSS1D} operations are $O(n_{\ell}^{1/2})$, including matrix application.
\item As indicated in Remark~\ref{rem:matrix-matrix}, for an HSS
matrix of size $k \times k$,  products of HSS and low-rank matrices
require $O(kq)$ work. Assumption~\ref{asm:compressed-blocks}  implies
all such products are $O(n_{\ell}^{1/2} \log_2(n_{\ell}))$. Low-rank
matrix matrix-vector multiplication has the same complexity.
\item  Finally,  both \textbf{LR\_to\_HSS1D} and \textbf{Rand\_ID}
involve interpolative decompositions of a matrix of size
$O(n_{\ell}^{1/2}) \times O(\log_2(n_{\ell}))$, and, therefore have complexity $O(n_{\ell}^{1/2} \log^2_2(n_{\ell}))$. Products between low-rank matrices are also of this complexity.
\end{enumerate}

\subsubsection*{Storage}
\begin{enumerate}
\item Again, since all HSS blocks behave as operators acting on one-dimensional 
box boundaries, storage is linear with respect to the number of nodes along the
boundary of the box: $O(n_{\ell}^{1/2})$.
\item A low-rank matrix of size $m \times n$ and rank $q$ occupies
$O((m+n)q)$ space in storage. By Assumption~\ref{asm:compressed-blocks}, storage of low rank blocks ($\ttA[L],\ttA[R]$ and off-diagonal blocks of F) is $O(n_{\ell}^{1/2} \log_2(n_{\ell}))$.
\end{enumerate}

Algorithms~\ref{alg:inter-lowrank}, \ref{alg:build-inv}, and \ref{alg:builde},
require only the operations listed above. We observe that all algorithms contain 
at least one $O(n_{\ell}^{1/2} \log^2_2(n_{\ell}))$ operation. In terms of storage, 
\textbf{INTER\_LOWRANK} and \textbf{BUILD\_Finv} store both HSS and low-rank blocks 
($O(n_{\ell}^{1/2} \log_2(n_{\ell}))$), and \textbf{BUILD\_E} one HSS block ($O(n_{\ell}^{1/2})$).
Hence, the \emph{compressed-block} interpolation operator build and inverse compression
algorithms perform $W^{CB}_{\ell} =
O(n_{\ell}^{1/2} \log^2_2(n_{\ell}))$ work per box and require $M^{CB}_{\ell} = O(n_{\ell}^{1/2} \log_2(n_{\ell}))$ storage for each set of matrices computed at level $\ell$. As a contrast, their \emph{dense-block} counterparts have $W^{DB}_{\ell} = O(n_{\ell}^{3/2})$ and $M_{\ell}^{DB} = O(n_{\ell})$.

We note that a more detailed complexity analysis may be performed to
obtain constants for each subroutine, given the necessary experimental
data about our kernel for a given accuracy $\acc$. The specific
dependance of these constants on accuracy is briefly discussed and
tested in Section~\ref{sec:vary-accuracy}.

We summarize the complexity estimates in the following proposition

\begin{proposition}
\label{prop:complexity}

{\textbf{Dense-Block Algorithms}} Let $\cA[A]$ be an $N \times N$ system matrix such that assumptions 1-3 hold.
Then, the dense-block tree build and inverse compression algorithms perform $O(N^{3/2})$ work. For NTI kernels, storage requirements and matrix apply are both $O(N \log N)$. For TI kernels, storage is $O(N)$, and matrix apply is  $O(N \log N)$.

\textbf{NTI Compressed-Block Algorithms} Let $\cA[A]$ be an $N \times N$, non translation invariant system matrix such 
that assumptions 1-4 hold. Then compressed-block tree build, inverse compression and HSS apply all perform $O(N)$ work, 
and require $O(N)$ storage.

\textbf{TI Compressed-Block Algorithms} Let $\cA[A]$ be an $N \times N$, 
translation invariant system matrix such that assumptions 1-4 hold. 
Then compressed-block tree build, inverse compression and HSS apply 
all perform $O(N)$ work, and require $O(N)$ storage. In fact, inverse compression work and storage are sublinear: 
$O(N^{1/2} \log^2_2 N)$ and $O(N^{1/2} \log_2 N)$, respectively.
\end{proposition}

The limitations of dense-block algorithms now become clear. With notation as in Lemma \ref{lemma:complexity},
dense-block algebra corresponds to $(p,q) = (3/2,0)$ for work and $(p,q) = (1,0)$ for storage, which precludes
overall linear complexity. For the compressed block algorithms, on the other hand, we have 
$(p,q)=(1/2,2)$ for work and $(p,q)=(1/2,1)$ for storage, which does yield linear complexity.

\begin{remark}
As we have previously observed, the compressed-block algorithm may be generalized to system matrices with other rank growth and behavior. More specifically, we observe that a sufficient condition to maintain optimal complexity could be that the outer HSS structure ranks grow as $O(n_{\ell}^p)$, for $p<1$, and that the most expensive operations such as the randomized interpolative decomposition also remain being sublinear. This would require all low-rank blocks to be of rank $q_{\ell} \sim n_{\ell}^{min\{1/3,(p-1)/2p\}}$.
\end{remark}

\begin{remark}
In all our practical implementations of the compressed-block algorithms, we perform dense computations for blocks up to a fixed threshhold skeleton size $k_{\rm cut}$, after which we switch to the fast routines. This may then be tuned as a parameter to further speed up these algorithms, and a straight-forward computation shows it does not alter their computational complexity.
\end{remark}

\section{Numerical Results}
\label{sec:numerical}
In this section, we describe a series of numerical experiments that test our inverse compression algorithm
on both translation-invariant (TI) and non-translation-invariant (NTI) kernels.

All results are obtained by running a Matlab implementation of our
algorithm are obtained on server nodes with Intel Xeon X5650, 2.67
GHz  processor, with no parallelization.

\subsubsection*{General formulation}

As we noted in the introduction, a considerable number of physical problems involve solving one or a system of integral equations of the form
\beqn
\cA[A] [\sigma](x) = a(x)\sigma(x) + \int_{\Omega} {b(x) \cA[K] (||x - y||) c(y)\sigma(y) dy} = f(x),
\eeqn
where $a$, $b$ and $c$ are given smooth functions, and $\cA[K](r)$ is related to a free-space Green's function. Then, given $\{ x_i \}_{i=1}^{N} \in \Omega = [-1,1]^2$ points on a regular grid with spacing $h$, we can perform a Nystr\"om discretization, obtaining a linear system $A\sigma = f$, with $A$ an $N \times N$ matrix with entries:

\beqn
A_{i,j} = a(x_i)\delta_{i,j} + h^2 b(x_i) \cA[K] (||x_i - x_j||) c(x_j)
\eeqn

\subsection{High accuracy: performance and scaling}
Most of our tests are done for a target accuracy of $\acc = 10^{-10}$ (referring to the local truncation error of all routines above); this is the ``stress test'' for our algorithm, as high accuracy requires
larger skeletons. We  measure the wall-clock timings and memory usage for the following parts of the solution process:

\begin{enumerate}
\item building the tree and interpolation operators;
\item inverse construction and compression;
\item inverse matrix-vector multiplication (solve, timings only).
\end{enumerate}

We compare against the dense-block HSS inversion algorithm.
We take $\cA[K](r)$ to be the 2D Laplace free-space Green's
function: $ \cA[K](r) = \frac{1}{2\pi} \log(r) $ and $a \equiv 1$.

\begin{itemize}
\item If $b \equiv c$, $A$ is symmetric, which leads to some computational savings in the HSS algorithms.
\item The case $b \not \equiv 1$ is our example of a \textbf{non translation invariant (NTI)} kernel. Such an equation appears in forward scattering problems (zero or low-frequency Lippman-Schwinger equation).
\item The case $b \equiv 1$ is one of our examples of \textbf{translation invariant (TI)} kernel. We recall that significant computational savings may be achieved in this case, given that we only build one set of matrix blocks per level.
\end{itemize}

We note that these results are typical for non-oscillatory Green's
function kernels: we have performed the same tests for the 3D Laplace
single layer potential and the 2D and 3D Yukawa Green's functions,
and found behavior analogous to the one described below.

\subsubsection{Non-translation-invariant kernel}

We take $b(x)$ to be a smooth function with moderate variation:
\beqn
b(x) = 1 + 0.5e^{-(x_1-0.3)^2-(x_2-0.6)^2}
\eeqn

Although this may seem like a simplistic choice, it does not matter much for computational performance
unless the problem is under-resolved or $b$ is $0$ on a large subdomain.

For the NTI case ranks of blocks in dense-block HSS matrices and ranks of low-rank blocks
vary both by level and inside each level. However, we observe that the rank growth for these blocks
matches our complexity assumptions.

We set leaf box size at $n_{max}=7^2$, and for problem sizes $N=784$ to $N=3211264$, we compare total inversion time (tree build + inverse compression) and total memory usage for the dense-block (HSS-D) and compressed-block (HSS-C) inverse compression algorithms.

\begin{table}[h!]
\centering
\begin{tabular}{ | r | c | c | c | c | }
\hline
  $N$ & HSS-D Time & HSS-C Time & HSS-D Memory & HSS-C Memory  \\
      & $O(N^{3/2})$ & $O(N)$ & $O(N \log N)$ & $O(N)$ \\
  \hline
  784 & 0.11 s  & 0.17 s  & 4.68 MB & 4.48 MB   \\
  3136 & 0.67 s & 1.70 s  & 29.09 MB & 25.24 MB  \\
  12544 & 4.50 s  & 8.32 s  & 159.59 MB & 123.07 MB   \\
  50176 & 31.45 s & 40.43 s  & 819.58 MB & 538.51 MB   \\
  200704 & 3.79 m  & 3.23 m  & 3.72 GB & 2.23 GB   \\
  802816 & 28.35 m  & 13.66 m  & 17.27 GB & 9.23 GB   \\
  3211264 & 3.58 hr  &  54.795 m  & 70.99 GB &  34.09 GB   \\
  \hline
\end{tabular}
\caption{Total inversion time and storage for the NTI inverse compression algorithms}
\label{tbl:NTI}
\end{table}

We observe a very close match between the experimental scaling and the
complexity estimates in Proposition \ref{prop:complexity}, which we recall on the second
row of Table~\ref{tbl:NTI}.

The slopes in a log-log plot (Figure \ref{fig:NTI_Inv}.A) show that
both inversion and tree build times are $O(N^{3/2})$ for the dense
block version and $O(N)$ for our accelerated algorithm. The break-even point
is around $N = 10^5$, for which both methods take about $1.5$
minutes. In the speedup plot
in \ref{fig:NTI_Inv}.B, we can clearly observe performance gains for
$N > 10^5$ due to difference in scaling. By $N = 10^7$ our algorithm gains one order of magnitude speedup.
Additional speed gains may be obtained by adjusting the size of the fine-scale boxes.

The slopes of the log-log plot (Figure \ref{fig:NTI_Inv}.C) confirm
that memory usage for these algorithms behaves like $O(N \log N)$ and
$O(N)$, respectively.  Our algorithm uses less memory in all
cases.(Figure~\ref{fig:NTI_Inv}.D).
By $N = 10^7$, it uses $2.5 \times$ less memory (150 GB, which amounts
to approximately 1800 doubles per degree of freedom).

\subsubsection*{Matrix compression} 
\begin{figure}[h!]
\begin{center}
\includegraphics[scale = 0.17]{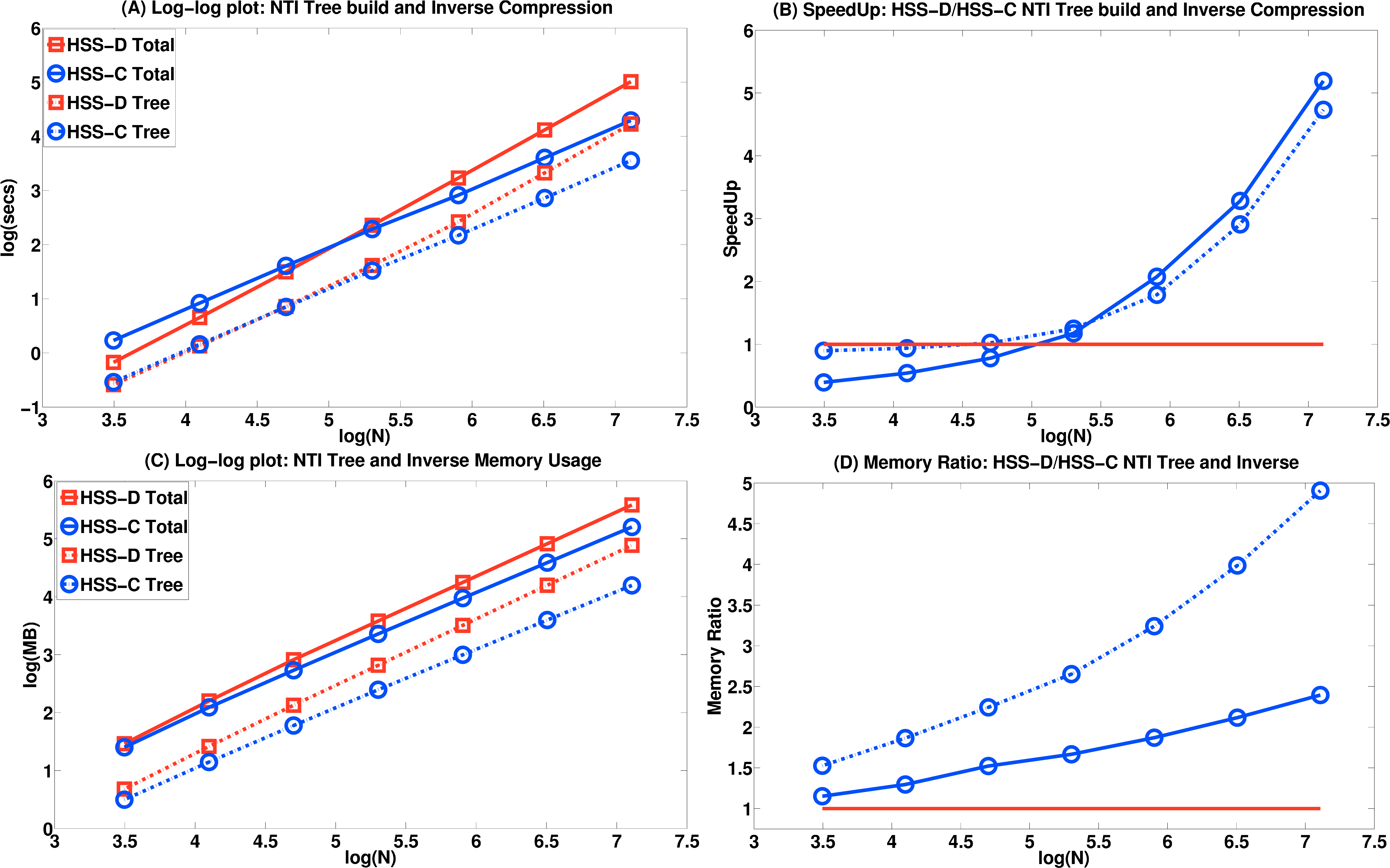}
\caption{\textit{NTI 2D Laplace Kernel: Inverse compression timings and memory usage. Total times and storage are represented in solid lines, Tree build and storage in dashed lines.}}
\label{fig:NTI_Inv}
\end{center}
\end{figure}

Figure \ref{fig:NTI_Inv}.A, shows that
the tree build time is about $18 \%$ of the total inversion
time. The constructed tree can also be used to compress the original
matrix $\cA[A]$ in $O(N \log N)$ work, and even in $O(N)$ work by
incurring the small additional cost of compressing sibling interaction
matrices as HSS 1D. Even though our focus on building a fast direct
solver, this suggests our accelerated method to compress interpolation
matrices $\ttA[L]$ and $\ttA[R]$ readily yields $O(N)$ matrix
compression. The compressed matrix can be used as a kernel-independent
FMM, with a significantly simplified algorithm structure.

\subsubsection{Translation-invariant kernel}
Exploiting translation invariance yields a substantial improvement in performance, which stems from the fact that only one set of matrices per level of $\cA[T]$ needs to be computed and stored.  The full-block HSS inverse compression still scales as $O(N^{3/2})$, but is $3\times$ faster for all tested problem sizes, and memory usage becomes more
efficient as $N$ grows. The lack of improvement in the asymptotic behavior of inverse compression when compared to the the non-translation-invariant case stems from the fact that $O(N^{3/2})$ work is performed on the top boxes of the tree in both cases.

As one would expect, our compressed-block algorithm displays sublinear
 scaling for both inverse compression time and storage, for the
 largest problems gaining an order of magnitude in both when compared
 to the NTI case.
 For $N = 10^7$, it has become about $11\times$ faster (29 minutes), using $30\times$ less memory (5 GB).
 Building the binary tree, which as mentioned amounts to
 compressing $\cA[A]$, displays similar performance gains. For $N = 10^7$, it takes only $4.5$ minutes to build and storage is only $1.3 GB$.

\begin{table}[h!]
\centering
\begin{tabular}{ | r | c | c | c | c | }
\hline
  N & HSS-D Time & HSS-C Time & HSS-D Memory & HSS-C Memory  \\
      & $O(N^{3/2}) $ & $O(N)$ & $O(N)$ & $O(N)$ \\
  \hline
  784 & 0.05 s  & 0.13 s  & 1.94 MB & 1.75 MB\\
  3136 & 0.21 s  & 0.98 s & 9.04 MB & 6.19 MB\\
  12544 & 1.40 s  & 3.41 s & 39.16 MB & 19.03 MB\\
  50176 & 9.68 s & 10.76 s & 163.19 MB & 52.09 MB\\
  200704 & 1.21 m  & 30.89 s & 666.39 MB & 151.41 MB\\
  802816 & 9.20 m  & 1.59 m  & 2.61 GB & 474.74 MB\\
  3211264 & 1.19 hr  &  6.68 m & 9.9359 GB &  1.56 GB   \\
  12845056 & 9.28 hr  &  29.22 m  & 39.74 GB &  5.29 GB \\
  \hline
\end{tabular}
\caption{Total inversion time and storage for the TI inverse compression algorithms}
\end{table}

\begin{figure}[h!]
\begin{center}
\includegraphics[scale = 0.17]{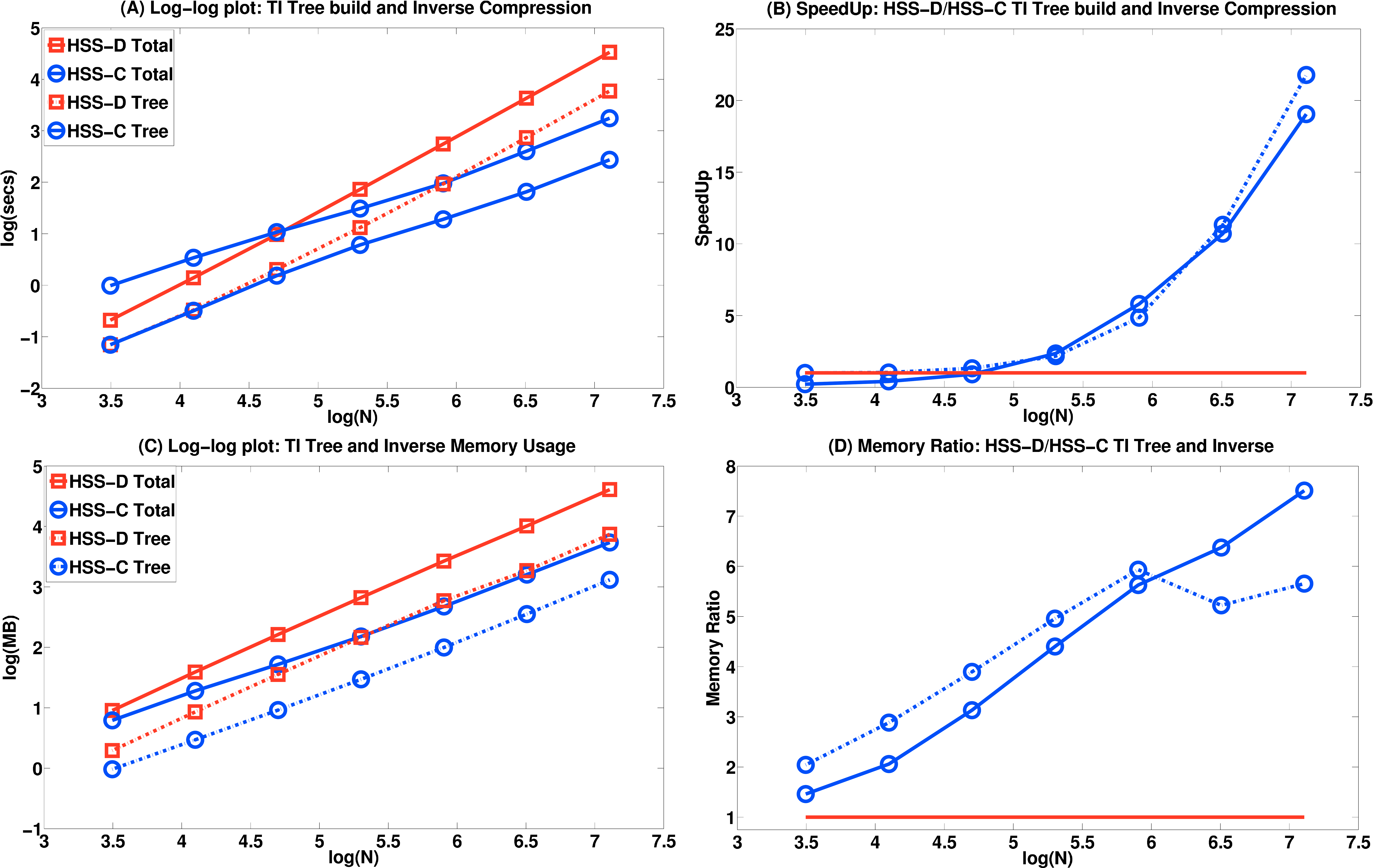}
\caption{\textit{TI 2D Laplace kernel: Inverse compression timings and memory usage. Total times and storage are represented in solid lines, Tree build and storage in dashed lines.}}
\label{fig:TI_Inv}
\end{center}
\end{figure}

The algorithm achieves parity with the dense-block HSS for lower
values of $N$ (around $N = 50000$), for which both methods take about
$10$ seconds to produce the HSS inverse. By $N =10^7$, our algorithm
is about $20 \times$ faster than the dense-block approach
(Figure \ref{fig:TI_Inv}). We again observe that tree build time is consistently about $15\%$-$20 \%$ of the total inversion time, and so our observation about matrix compression holds.
Figure \ref{fig:TI_Inv}.C shows sublinear scaling for memory usage. It
remains true that our algorithm uses less memory in all cases. For $N
=10^7$, it uses $8 \times$ less memory
(5 GB, which amounts to 50 doubles per degree of freedom).

\subsubsection{Inverse matrix-vector multiplication}
In order to test the inverse matvec algorithm, we run it with  ten
random right-hand sides and compute time per solve (in seconds) for
both NTI and TI examples. We note that, given that the same matrix
blocks are applied to the same subsets of each right-hand side, code
can be easily vectorized leading to performance gains if multiple
solves are to be performed simultaneously. If our matrix is
translation-invariant additional speedup result from the fact that the
same set of matrices is applied to vectors corresponding to all boxes in a given level.

We note that for all cases, the inverse apply is very fast, scaling
linearly and remaining well under a minute for sizes up to $N \sim
10^7$.  As expected, the differences between the standard and compressed block apply are small, the latter becoming marginally faster around $N \sim 3\times10^6$.

\begin{table}[h!]
\centering
\begin{tabular}{ | r | c | c | c | c | c | c |}
\hline
  N & NTI HSS-D & NTI HSS-C  & TI HSS-D & TI HSS-C  &  NTI HSS-C & TI HSS-C  \\
  & $O(N \log N) $ & $O(N)$ & $O(N \log N)$ & $O(N)$ & Error & Error \\
  \hline
  784 & 0.0014  & 0.0018  & 0.0007 & 0.0011 & 1.6e-14 & 9.2e-15  \\
  3136 & 0.0064  & 0.0090  & 0.0031 & 0.0046 & 1.8e-14 & 1.8e-14 \\
  12544 & 0.0292  & 0.0362  & 0.0137 & 0.0162 & 8.6e-11 & 5.7e-11 \\
  50176 & 0.1320  & 0.1546  & 0.0590 & 0.0600 & 1.6e-10 & 1.7e-10  \\
  200704 & 0.5993  & 0.6772  & 0.2819 & 0.2512 & 2.3e-10 & 1.6e-10 \\
  802816 & 2.6611  & 2.8193  & 1.2709 & 1.0763 & 4.0e-10 & 3.8e-10  \\
  3211264 & 11.816  &  11.737  & 5.77296 &  4.5650 & 5.1e-9 & 1.6e-9 \\
  12845056 & \color{red}52.468  &  \color{red}48.8641  & 25.8312 &  19.3619 & - & 3.2e-9  \\
  \hline
\end{tabular}
\caption{Inverse apply timings (in seconds) for both NTI and TI algorithms. Numbers in red are extrapolated from previous data.}
\end{table}

For all cases, if $\mathcal{A}$ and $\mathcal{A}^{-1}$ are our compressed, approximate matrix and inverse, we use the following error measure:
\beqn
E = \frac{|| v - \mathcal{A}\mathcal{A}^{-1}v||}{||v||}
\eeqn

Taking the maximum over a number of randomly generated right-hand-sides $v$. In this measure, which can be thought of as an approximate residual, both algorithms achieve the desired target accuracy for the inverse apply.
We may then bound the exact residual by:

\beqn
|| v - A\mathcal{A}^{-1}v || \leq E + || \mathcal{A} - A || || \mathcal{A}^{-1} v||
\eeqn

\subsection{The effect of varying accuracy}
\label{sec:vary-accuracy}

If we set a lower target accuracy, this can significantly speed up both algorithms presented above:  the number of boundary layers needed,  the size of skeleton sets in $\cA[T]$, ranks in low-rank and HSS 1D blocks in the compressed-block algorithm are all lower, and low-rank factorization can be performed faster.

In this set of tests, using the Laplace single layer potential in 2D, we compare their behavior for target accuracies $\acc = 10^{-5}$ and $10^{-10}$.

\subsubsection*{Dense-block algorithm} In this case, only the change in skeleton set size is relevant. As in Section~\ref{sec:complexity}, we can bound $k_i \le k_{\ell} = C_{\acc}n_{\ell}^{1/2}$. In particular,
$\acc=10^{-5}$ requires one boundary layer and $\acc=10^{-10}$ requires
two, and so $C_{10^{-10}} \sim 2C_{10^{-5}}$, as expected from standard multipole estimates ($C_{\acc} \sim \log(1/\acc)$).
Following our complexity analysis, we observe both tree build and inverse compression algorithms perform $O(k_i^3)$ work per box, and so the constant in the leading term is of the form $\hat{C}\log^3(1/\acc)$. This would imply a factor of $8$ between $\acc=10^{-5}$ and $\acc=10^{-10}$.
Analogously, inverse matvec work and memory usage are $O(k_i^2)$ per box, and so we expect a factor of 4 between these two cases.

\subsubsection*{Compressed-block algorithm}
Performance of the compressed block algorithm depends on all  three quantities mentioned above:
the size of 2D HSS skeletons is again $O(\log(n_{\ell}/\acc))$, but
fast operations are performed on these.
Hence, we expect a factor of $\log(1/\acc)$ to appear, at worst.
Assuming sizes of dense blocks in low-rank and HSS 1D representations
are $O(\log(n_{\ell}/\acc))$, the most expensive one-dimensional
dense-block HSS perations are linear with a constant $O(\log^3(1/\acc))$, and again matvecs and storage get a factor  $O(\log^2(1/\acc))$.

Hence, we expect the factors to behave as $O(\log^4(1/\acc))$ ($16 \times$) for tree build and inverse compression, and as $O(\log^3(1/\acc))$ ($8 \times$) for inverse apply and memory storage.

\begin{table}[h!]
\centering
\begin{tabular}{ | r | c | c | c | c | c | c | c | c |}
\hline
  N & HSS-D Time & Ratio & HSS-C Time & Ratio & HSS-D Memory & Ratio & HSS-C Memory & Ratio \\
  \hline
  784 & 0.01 s & 3.1 & 0.03 s  & 6.7 & 0.45 MB & 3.3 & 0.44 MB & 4.1 \\
  3136 & 0.03 s  & 5.2 & 0.07 s & 15 & 1.92 MB & 3.7 & 1.75 MB & 4.2 \\
  12544 & 0.20 s  & 6.4 & 0.75 s & 6 & 7.98 MB & 3.8 & 5.40 MB & 4.3 \\
  50176 & 1.35 s & 6.9 & 2.83 s & 6 & 32.53 MB &3.9 & 13.67 MB & 4.4 \\
  200704 & 9.77 m  & 7.4 & 7.87 s & 6.4 & 131.39 MB & 4 & 36.59 MB & 4.4 \\
  802816 & 1.22 m  & 7.6 & 20.88 s  & 7.3 & 528.19 GB & 4 & 107.38 MB & 4.4 \\
  3211264 & 9.13 m  & 7.9 &  1.19 m & 8 & 2.07 GB & 4 &  349.50 MB & 4.2 \\
  12845056 & 1.14 hr  & 8.1 &  4.56 m  & 8.9 & 8.34 GB & 4 &  1.21 GB & 3.8 \\
  \hline
\end{tabular}
\caption{ Inverse compression time and storage of the TI algorithms
  for $\acc = 10^{-5}$, and the ratio of the quantities corresponding to $10^{-10}$ and $10^{-5}$ target accuracy}
\end{table}

Experimental results show that our estimates for the dense block
algorithm (HSS-D) are close: as N grows, we observe that the ratio
between compression times converges to $\mathbf{8}$, and that for
memory usage to $\mathbf{4}$, as expected. However, although there  more variation in the compressed-block results, we generally observe our estimates to be conservative, since the ratios are fairly similar to the dense-block case.

Our hypotheses also appear to be conservative for the inverse apply algorithms: ratios between these two accuracies tend to be between $\mathbf{2}$ and $\mathbf{3}$ for all cases tested.

We omit the results for the non-translation-invariant case, since the comparison and ratios between both target accuracies are quite similar. Overall, this implies that the analysis and observations in Section 5.1 can be applied with little modification to the case $\acc =  10^{-5}$: for both algorithms, a faster but less accurate inverse can be compressed $\mathbf{8}$ times faster, stored using $\mathbf{4}$ times less memory, and can be applied two or three times faster than the high accuracy case.

\subsection{Low Frequency Oscillatory Kernels}

Finally, we perform the same set of experiments as in Section 5.1, but this time taking $\mathcal{K}(r)$ in our general formulation to be the 2D Helmholtz free-space Green's function for wave number $k$: 

\beq
\mathcal{G}_{k}(r) = -\frac{i}{4} H^1_0(k r)
\label{eqn:H2D_kernel}
\eeq

where $H^1_0$ is the Hankel function of the first kind. As we will see in Section 5.4, this is akin to solving the Lippmann-Schwinger equation for a given frequency $k$. 

As mentioned before, with some minor modifications our solver is able to handle oscillatory kernels for low frequencies. We briefly note the main differences with non-oscillatory problems:

\begin{itemize}
\item{\textbf{Skeleton Structure:}} For high accuracy, skeletons with
consisting of more layers of points are needed in comparison with non-oscillatory kernels such as Laplace. Besides this, as the wavenumber grows it is also necessary to keep points inside of the box in the skeleton set. Experimentally we determine that keeping several points per wavelength is enough to maintain accuracy for matrix and inverse compression.

\begin{figure}[h!]
\begin{center}
\includegraphics[scale = 0.09]{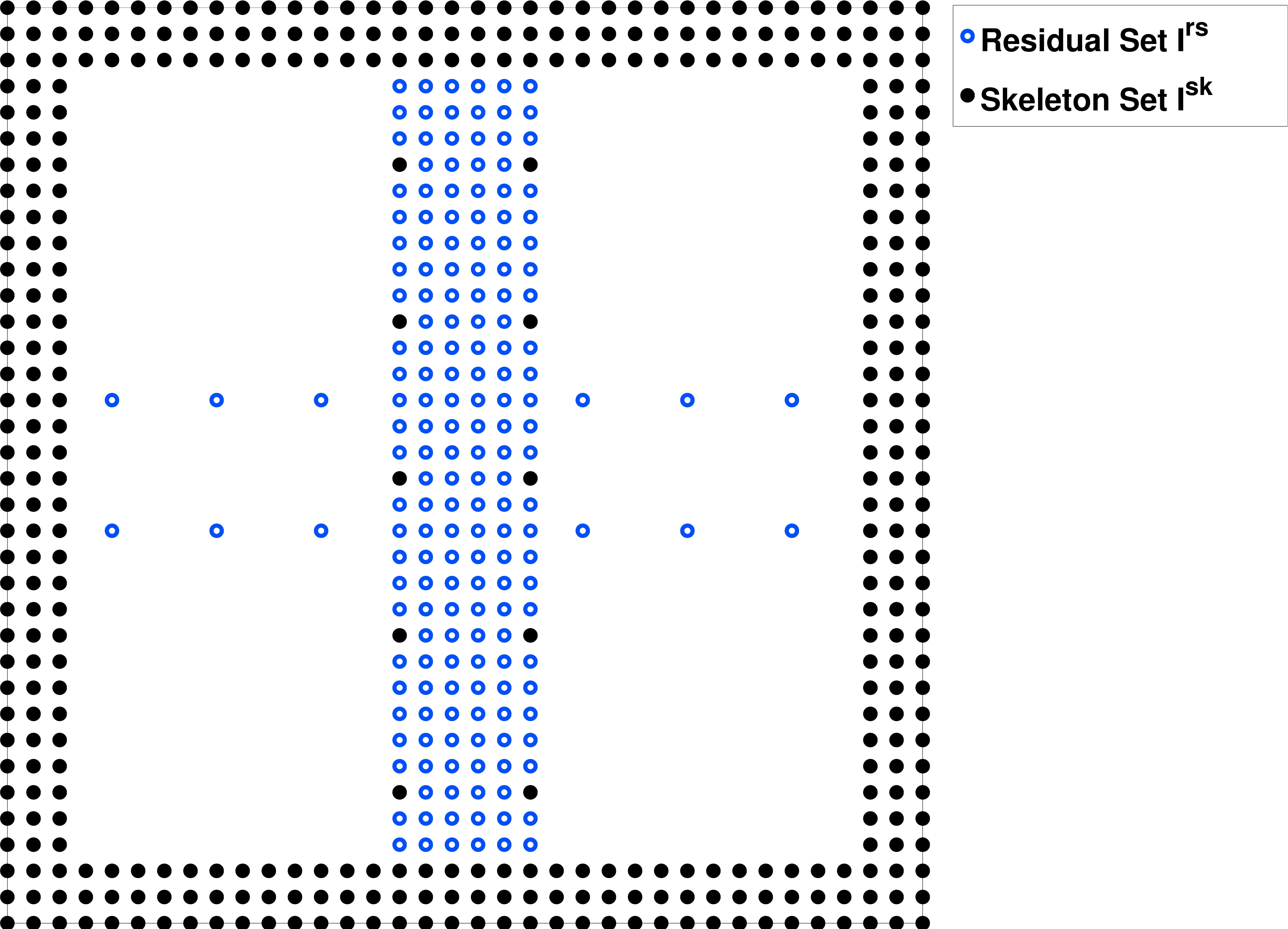}
\caption{\textit{2D Helmholtz kernel skeleton set structure: Skeleton points are in black and residual points in blue. A few points per wavelength are kept on the interface between the two children boxes.}}
\label{fig:H2D_skeletons}
\end{center}
\end{figure}

\item{\textbf{Conditioning:}} The Lippmann-Schwinger equation is moderately
ill-conditioned, with condition number depending on $k$. This impacts
the matrix inversions in both levels of compression resulting
several digits of accuracy.  We note that this is mild for the
wavenumbers tested, and that we can always gain them back by one round of iterative refinement (matrix and inverse applies remain quite fast).
\item{\textbf{Rank growth assumptions:}} Ultimately, as $k$ grows, the assumptions stated in Section 4 will fail: The size of skeletons in compressed blocks will depend on $k$, as well as the number of points needed in our discretization.
\end{itemize}

This points to the fact that a different approach is needed to build a fast direct solver that can handle moderate and high frequency regimes.

\subsubsection{Inverse compression results: TI and NTI cases}

We present experimental results for both the non-translation-invariant and the translation-invariant cases. Let $\kappa = k / 2 \pi$ be the number of wavelengths that fit in our domain, the unit box. For $\kappa = 4,8,16,32$, we set a leaf box size at $n_{max} = 9^2$ and test for increasing problem size $N$. Here we compare total inversion time and memory usage for the HSS-D and HSS-C inverse compression algorithms.

We observe that, while the cases for $\kappa=4,8$ behave similarly to their Laplace counterpart (skeletons consist of 2 boundary layers), for $\kappa=16,32$ and higher at least 3 layers are needed, and accuracy in matrix compression deteriorates unless we also keep a few extra points per wavelength inside the box. As we will see, this change in skeleton set structure is the main cause for differences in behavior between these two sets of examples.

We note that, given that the Helmholtz kernel is complex-valued, two doubles are stored for each matrix entry. This implies that memory usage can be a priori expected to be at least 2 times bigger than that of real-valued kernels.

\begin{figure}[h!]
\begin{center}
\includegraphics[scale = 0.15]{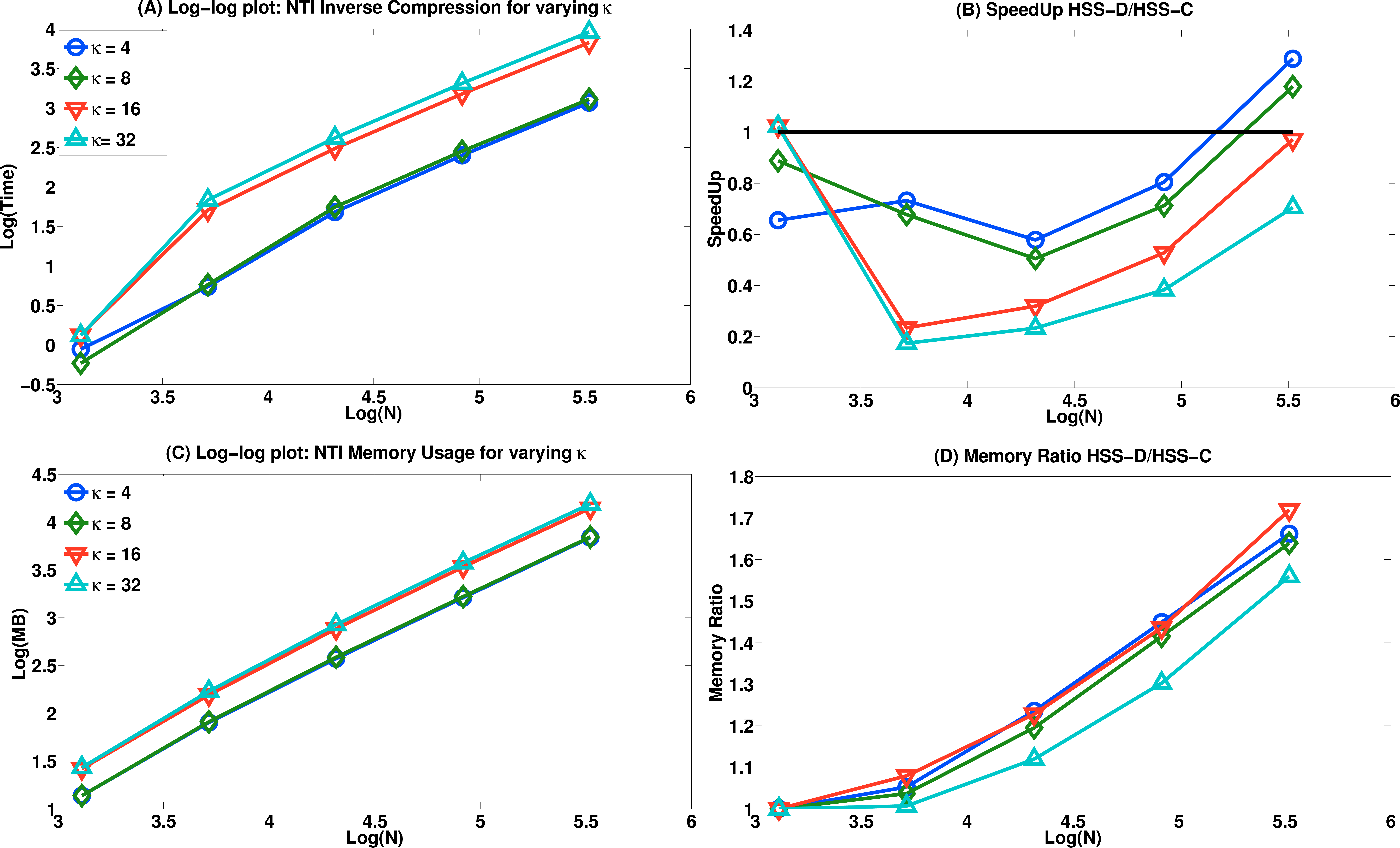}
\caption{\textit{NTI 2D Helmholtz kernel: Inverse compression timings and memory usage for increasing values of $\kappa$}}
\label{fig:H2D_NTI_Inv}
\end{center}
\end{figure}

We first observe that for all cases, experimental scaling again coincides with the expected complexity (the $O(N^{3/2})$ HSS-D algorithm is not plotted in figure \ref{fig:H2D_NTI_Inv} to avoid clutter). For the HSS-C algorithm, slopes in \ref{fig:H2D_NTI_Inv}.A quickly approach to 1 as N grows.

\textbf{Break-even points and speedup:} on the speedup plot in figure \ref{fig:H2D_NTI_Inv}.B, we can observe that the point where the compressed-block algorithm overcomes its dense-block counterpart grows slightly with $\kappa$. While it is around $10^5$ for $\kappa = 4,8$, it grows closer to $10^6$ for $\kappa=16,32$. We can also see a moderate speedup is gained after these break-even points due to better scaling, similarly to the case for Laplace.

\textbf{Memory Usage:} In figures \ref{fig:H2D_NTI_Inv}.C and \ref{fig:H2D_NTI_Inv}.D  we again can observe the HSS-C algorithm always provides extra compression in terms of storage, and this improves as N grows. By $N \sim 10^6$, it uses $3-4$ times less memory than HSS-D.

For $N \sim 10^6$ and $\kappa = 4,8$, the HSS-C algorithm takes about $1.3$ hours to produce the HSS inverse, requiring $27GB$ of storage. For $\kappa = 16,32$, it takes 7 and 10 hours to build the inverse, respectively, and storage requirements go up to $\sim 60GB$. Comparing its performance with the Laplace NTI kernel in Section 5.1 (which requires 0.5 hours and 10 GB for this size), we see that as the wave number increases, it takes considerably more time and storage to compress the inverse. The most dramatic change can be observed for $\kappa = 16,32$, since skeleton sets become much bigger in size.

\begin{figure}[h!]
\begin{center}
\includegraphics[scale = 0.15]{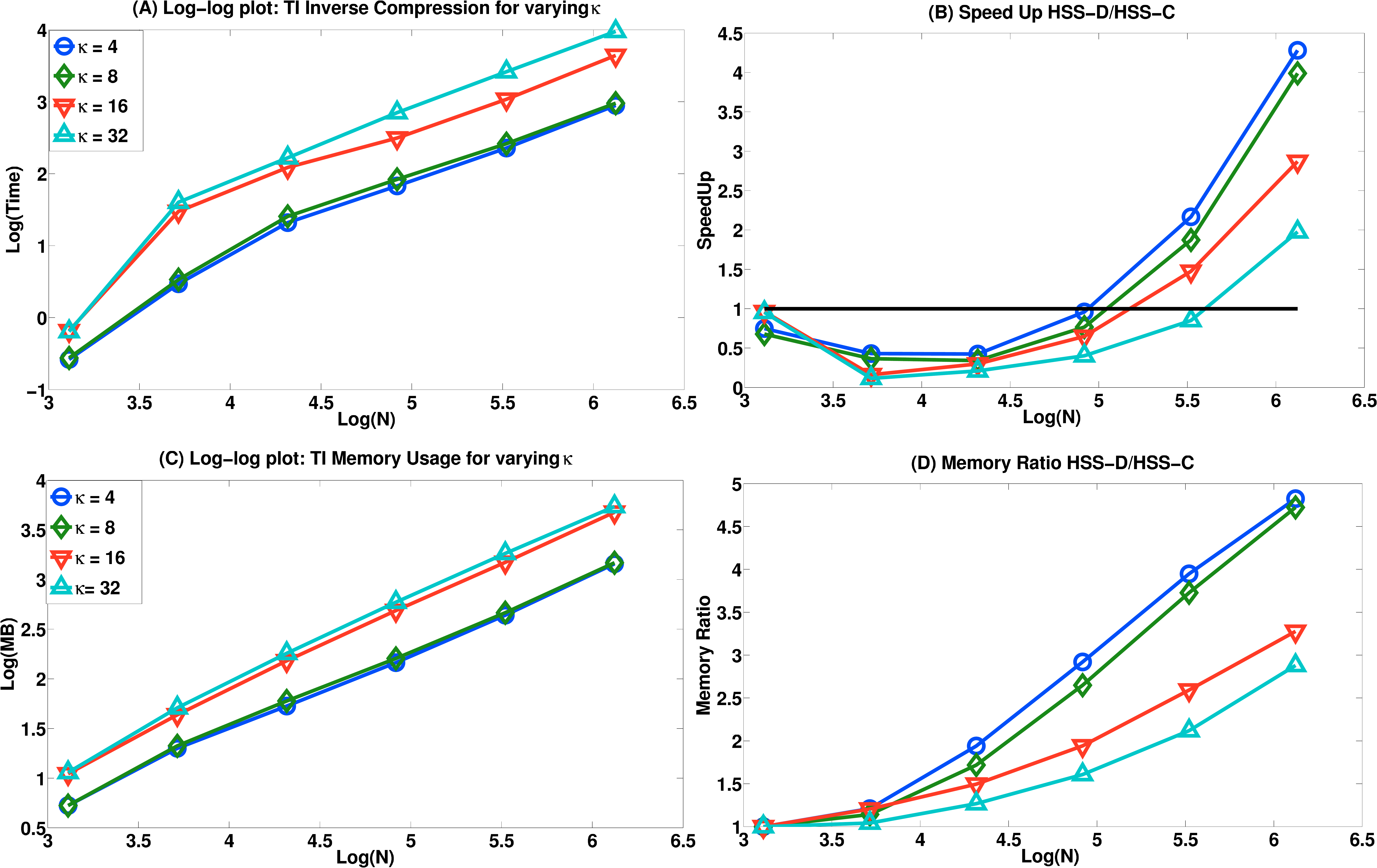}
\caption{\textit{TI 2D Helmholtz kernel: Inverse compression timings and memory usage for increasing values of $\kappa$}}
\label{fig:H2D_TI_Inv}
\end{center}
\end{figure}

Exploiting translation invariance yields significant performance gains for all tested wavenumbers. The HSS-D algorithm is again $O(N^{3/2})$, and about $3 \times$ faster than the NTI case. We also observe sublinear scaling for inverse compression and storage for the HSS-C version, although a bit less dramatic than in the case for the Laplace kernel. By $N \sim 10^6$, it has become $5-6 \times$ faster, using $15 \times$ less memory.

Since only one set of matrices is computed for each level in $\mathcal{T}$, we can observe more clearly the impact of adding an extra layer of skeleton points: inverse compression time and storage become about 1 order of magnitude higher. Also, there is a more pronounced difference between $\kappa=16$ and $\kappa=32$ for bigger problem sizes.

\textbf{Break-even points and speedup:} On the speedup plot in
figure \ref{fig:H2D_TI_Inv}.B, we see that break-even points are a bit
smaller when compared with the NTI case, and that they again grow with
$\kappa$. For $\kappa=4$ it is below $10^5$, for $\kappa = 32$, it
happens around $3\times 10^5$. Speedups are again considerably better and faster growing due to sublinear scaling, especially for small wavenumbers.

\textbf{Memory Usage:} The gain in compression is more rapid, especially for low wavenumbers. By $N \sim 10^6$, it uses $3-5$ times less memory than HSS-D. Also, in \ref{fig:H2D_TI_Inv}.D, we see the slopes for the memory ratio are higher than in the NTI case.

For $N \sim 10^6$ and $\kappa = 4,8$, the HSS-C algorithm takes about $15$ minutes to produce the HSS inverse, requiring $1.5GB$ of storage. For $\kappa = 16,32$, these go up to 1.2 and 2.6 hours to build the inverse, taking $4.6$ and $5 GB$ of storage . Comparing its performance with the Laplace TI kernel (2 minutes and 0.5 GB for this size), we observe similar differences in storage as in the NTI case, but much more drastic ones in terms of inverse compression time.

\subsubsection{Inverse matrix-vector multiplication}
Finally, we test the inverse apply algorithm as in Section 5.1.3. Since the differences between dense and compressed block are quite small, we present results for the HSS-C apply only, for $\kappa=8$ and $\kappa=32$. The inverse apply is still considerably fast and retains optimal scaling.

\begin{table}[h!]
\centering
\begin{tabular}{ | r | c | c | c | c | c | c |}
\hline
  N & NTI HSS-C & NTI HSS-C  & TI HSS-C & TI HSS-C & TI HSS-C & TI HSS-C \\
  & $\kappa=8$ & $\kappa=32$ & $\kappa=8$ & $\kappa=32$ & Error ($\kappa=8$) & Error ($\kappa=32$)\\
  \hline
  1296 & 0.0036  & 0.0074 & 0.0022 & 0.0048 & 1e-13 & 1e-12 \\
  5184 & 0.0205  & 0.0596  & 0.0126 & 0.0365 & 5.5e-10 & 2.6e-8 \\
  20736 & 0.1102  & 0.3046  & 0.0587 & 0.1663 & 1.6e-10 & 8.7e-9 \\
  82944 & 0.5197  & 1.4962  & 0.2465 & 0.8148 & 2.3e-9 & 3.6e-8 \\
  331776 & 2.5431  & 6.7812  & 1.0350 & 3.4867 & 1.0e-10 & 2.0e-8 \\
  1327104 & \color{red}10.1684  & \color{red}27.1248  & 4.3536 & 15.4991 & 1.2e-9 & 4.9e-8 \\
  \hline
\end{tabular}
\caption{Inverse apply timings (in seconds) for both NTI and TI algorithms, for different values of $\kappa$. Numbers in red are extrapolated from previous data.}
\end{table}

For $\kappa=16,32$, we usually lose two digits in our accuracy measure ($ \sim 10^{-8}$), when compared to our target $10^{-10}$. As was mentioned above, this is an effect of conditioning and can be addressed if necessary by cranking up accuracy or by one round of iterative refinement.

\subsection{2D scattering problem: Lippmann-Schwinger equation}

\subsubsection*{Background on scattering} Consider an acoustic scattering problem in $\mathbb{R}^2$ involving a ``soft''
scatterer contained in a domain $\Omega$. For simplicity, we assume that the ``incoming field'' is generated from a point
source at the points $x_s \in \Omega^{\rm c}$. A typical mathematical model would take the form:
\beq
-\Delta u(x) - \frac{\omega^2}{v(x)^2} u(x) = \delta (x-x_s),\qquad  x \in \mathbb{R}^2
\label{eqn:helmholtz_freespace}
\eeq
where $\kappa$ is the frequency of the incoming wave, and where $v(x)$ is the wave-speed at $x$.
We assume that the wave-speed is constant outside of $\Omega$, so that $v(x) = v_0$ for
$x \in \Omega^{\rm c}$. To make the equation well-posed, we assume that $u(x)$ satisfies
the natural ``radiation condition'' at infinity. We define the ``wave number'' $k$ via
\beqn
k = \frac{\omega}{v_0},
\eeqn
and a function $b = b(x)$ that quantifies the ``deviation'' in the wave-number from the free-space wave-number via
\beqn
b(x) = k^2 - \frac{\omega^2}{v(x)^2}.
\eeqn
Observe that $b(x) = 0$ outside $\Omega$. Equation (\ref{eqn:helmholtz_freespace}) then takes the form
\beq
-\Delta u(x) - k^2 u(x) + b(x)\,u(x) = \delta (x-x_s),\qquad  x \in \mathbb{R}^2.
\label{eqn:scatter_local}
\eeq
Finally, let $G_k$ denote the Helmholtz fundamental solution as in equation \ref{eqn:H2D_kernel},
and split the solution $u$ of (\ref{eqn:scatter_local}) into ``incoming'' and ``outgoing'' fields
so that $u = u_{in} + u_{out}$, where
\beqn
u_{in}(x) = G_k(x,x_s).
\eeqn
Then $u_{out}$ satisfies:
\beq
-\Delta u_{out}(x) - k^2 u_{out}(x) + b(x)\,u_{out}(x) = -b(x)\,u_{in}(x),\qquad  x \in \Omega.
\label{eqn:outgoing}
\eeq

When converting (\ref{eqn:outgoing}) as an integral equation, we assume that $b(x)$ is non-negative,
and look for a solution of the form
\beq
u_{out}(x) = [S_k (\sqrt{b}\,\tau)](x) = \int_{\Omega} {G_k(x,y) \sqrt{b(y)}\,\tau(y)\,dy}.
\label{eqn:ansatz}
\eeq
Inserting (\ref{eqn:ansatz}) into (\ref{eqn:outgoing}), and dividing the result by $\sqrt{b(x)}$, we find the equation
\beq
\tau(x) + \int_{\Omega} {\sqrt{b(x)}\,G_k(x,y)\,\sqrt{b(y)} \tau(y) dy} = -\sqrt{b(x)}\, u_{in}(x),\qquad x \in \Omega
\label{eqn:LS}
\eeq

\subsubsection*{Test problem} For our numerical tests, we take:
\beqn
b(x_1,x_2) = 0.25 k^2 (1 + \tanh D(1 - \epsilon -|x_1|)) (1 + \tanh D(1 - \epsilon -|x_2|))
\eeqn
Note that
$b$ is equal to $1$ in a subdomain of $[-1,1]^2$ and transitions exponentially to $0$ close to the boundary of the square. Parameters $D$ and $\epsilon$ can be manipulated to change this layer's position and width. We note that if $b$ close to $0$ up to machine precision, the tree build routine must be modified to pick skeleton points along the support of $b$.

\begin{figure}[H]
\begin{center}
\includegraphics[scale = 0.05]{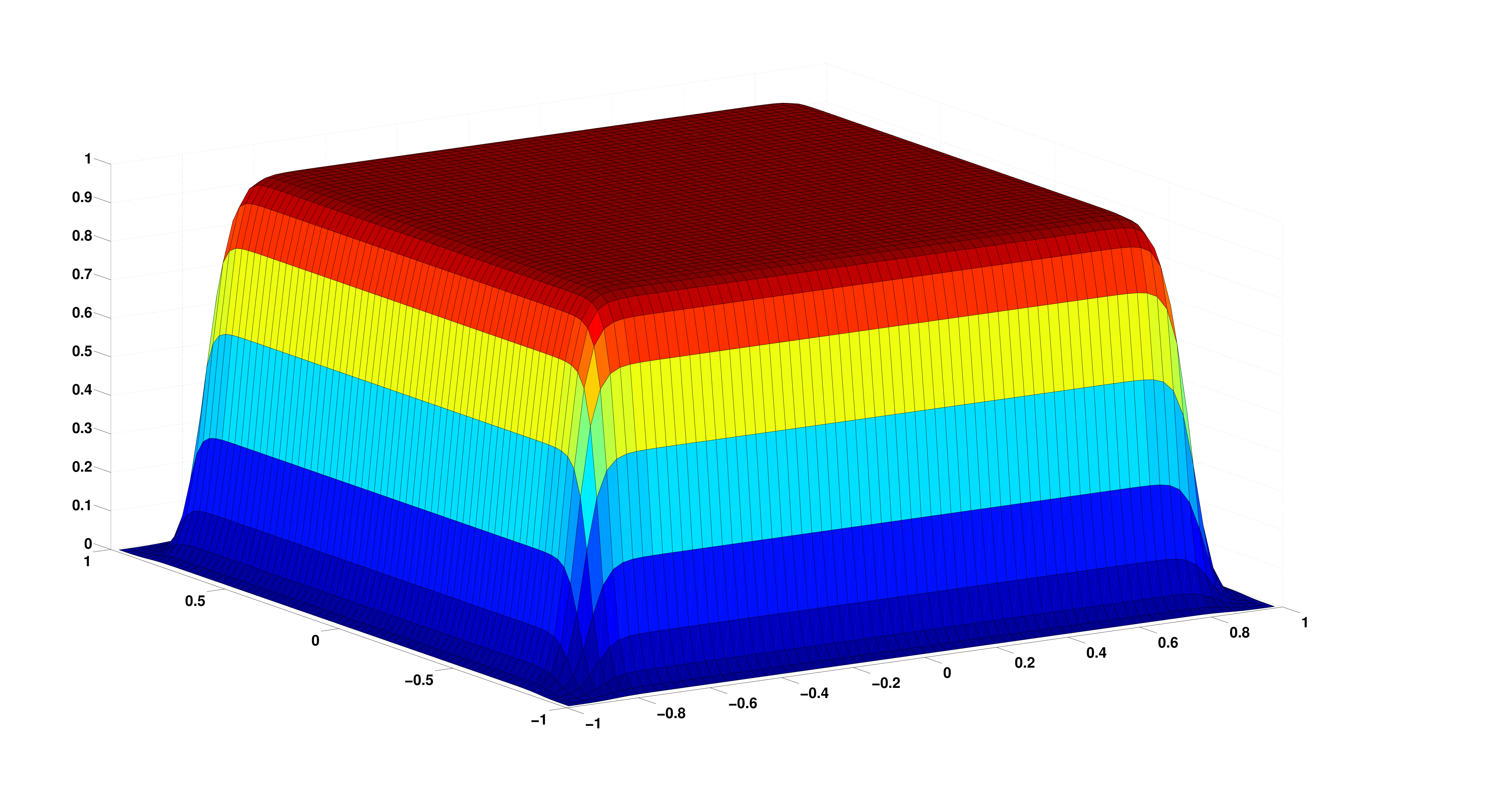}
\caption{\textit{Plot of deviation $b$ for the scatterer in our Lippman-Schwinger equation}}
\label{fig:LS_bump}
\end{center}
\end{figure}

We implement an $O(h^4)$ accurate corrected trapezoidal rule \cite{duan2009high} by adding a diagonal correction to our kernel. Since our algorithm compresses off-diagonal blocks, this has no bearing on the performance of our direct solver. We note that higher order (up to $O(h^{10})$) quadratures of this form can be applied if needed.

\subsubsection*{Numerical results} For $\kappa = 4,8$ and increasing problem size $N$, we solve the Lippmann-Schwinger equation 
(\ref{eqn:LS}) using our solver for non-translation invariant, symmetric operators.

For both for the approximate matrix $\mathcal{A}$ and its inverse, we first measure the empirical order of convergence when applying them to the same right hand side while refining grid size $h$. We confirm that the order is approximately $O(h^4) = O(N^{-2})$ until the error is comparable to the target accuracy.

Since there is a much higher contrast in the kernel entries, and they get close to $0$ at rows and columns corresponding to the box boundary, through preliminary experiments we observe that $2$ layers yield moderate accuracy ($\sim 10^{-6}$), and $3$ are needed for high accuracy ($\sim 10^{-10}$). We note that for cases where there is small variation within the domain of interest $\Omega$, this approach could be optimized by giving special treatment to boxes which intersect its boundary.

Also, we observe that the inverse of the translation-invariant operator in Section 5.3 may be used as a right preconditioner for this problem. Using BiCGstab, we find that for a given wave number $\kappa$, the number of iterations to solve the preconditioned system is moderate and independent of problem size $N$.

We test the direct solver and the preconditioned iterative solver proposed above (Table \ref{tbl:LS_inv}) For the latter, the compression step includes tree-build (matrix compression) for the Lippmann-Schwinger kernel, and inverse compression for the corresponding translation-invariant Helmholtz kernel (our preconditioner). 

\begin{table}[h!]
\centering
\begin{tabular}{ | r | c | c | c | c | c | c |}
\hline
  N & HSS-D Time & HSS-C Time & Iter Time & HSS-D Memory & HSS-C Memory & Iter Memory \\
  \hline
  784       & 1.92 s     & 1.86 s     & 1.72 s    & 18.74 MB    & 18.74 MB  & 9.63 MB     \\
  3136     & 11.81 s   & 24.93 s   & 9.30 s    & 125.87 MB    & 92.30 MB  & 54.37 MB    \\
  12544   & 1.31 m    & 2.68 m   & 1.13 m   & 720.90 MB & 464.81 MB & 258.97 MB  \\
  50176   & 8.03 m    & 14.20 m & 5.65 m   & 3.69 GB      & 5.56 GB    & 1.19 GB       \\
  200704 & 47.13 m  & 1.17 h    & 25.71 m & 18.34 GB    & 13.22 GB    & 5.53 GB       \\
  802816 & 4.78 h     & 5.70 h    & 1.93 h     & 90.36 GB   & 47.16 GB   & 25.38 GB     \\
  \hline
\end{tabular}
\caption{ Total inverse compression time and storage for the Lippmann-Schwinger equation
  for $\kappa = 8$ and $\acc = 10^{-10}$ target accuracy}
\label{tbl:LS_inv}
\end{table}

Finally, on Table \ref{tbl:LS_apply} we compare inverse apply (solve) stage for both algorithms and for the iterative approach.

\begin{table}[h!]
\centering
\begin{tabular}{ | r | c | c | c | c | c | c |}
\hline
  N & HSS-D Solve & HSS-C Solve & HSS Error & Iter Solve & Iter $\#$ & Error \\
  \hline
  784       & 0.03 s   & 0.02 s    & 2.7e-11  & 0.59   s    & 20      & 8.0e-9  \\
  3136     & 0.08 s   & 0.17 s    & 5.2e-9    & 3.13   s    & 30      & 9.2e-9  \\
  12544   & 0.36 s   & 0.82 s    & 7.1e-8    & 14.34 s    & 28      & 4.4e-9  \\
  50176   & 1.74 s   & 3.73 s    & 3.6e-11  & 56.12 s    & 28      & 3.6e-9   \\
  200704 & 9.54 s   & 20.21 s  & 9.3e-11  & 236.23 s  & 29.5   & 7.7e-9   \\
  802816 & 52.50 s & 109.32 s & 1.2e-10 & 1010.1 s  & 29.5   & 4.4e-9  \\
  \hline
\end{tabular}
\caption{Inverse apply timings for the Lippmann-Schwinger equation for $\kappa = 8$ and $\acc = 10^{-10}$ target accuracy}
\label{tbl:LS_apply}
\end{table}

We note that for the HSS-C algorithm, one round of adaptive refinement is needed in order to attain the target accuracy for the approximate residual. This means two inverse and one matrix applies are used in the solve.

\section{Conclusions and Future Work}
\label{sec:conclusions}
We have described a direct solver for volume integral equations in the plane that attains optimal $O(N)$ complexity and high practical efficiency for problems with non-oscillatory (or moderately oscillatory) kernels. 
The solver displays high performance for large problem sizes, even for high target accuracies. It gains a significant additional advantage when dealing with translation invariant kernels, since only one set of structured matrices needs to be computed for each level of the tree, resulting in sublinear scaling of inverse compression time and storage.

The solver is based on the recursive skeletonization scheme described in \cite{MR2005}, which would have $O(N^{3/2})$ complexity if applied to the volume integral equations considered here. We attained the acceleration to optimal $O(N)$ complexity by using structured matrix algebra to manipulate certain large dense matrices that arise in the computation. Specifically, we used the so called \emph{Hierarchically Semi-Separable (HSS)}, format which exploits rank deficiencies in the off-diagonal blocks of the matrix. As noted in Section 5, the accelerated direct solver provides better compression than the dense-block algorithm of \cite{MR2005} in all tested environments, and it becomes faster for sizes greater than $N \sim10^5$. Because of optimal scaling, this implies a one order of magnitude speedup for every increase in two orders in problem size $N$.

We also note that the increase in memory efficiency of our algorithm is of significant importance: given that all structured-matrix and FMM-like methods are memory intensive, the lack of sufficient compression can render them prohibitive for large enough problems.

The fast tree build routine described can be modified into an $O(N)$ matrix compression algorithm which is competitive with the FMM and may be used with iterative methods such as GMRES, see Section 5.

\paragraph{Future work directions.} There are several directions we 
plan to explore: 
\begin{itemize}
\item{\textbf{Extension to surfaces in} $\mathbf{3D}$:} we expect that the compressed-block HSS algorithms presented in this paper can be ported with no major structural differences to the case of solving boundary integral equations on surfaces in 3D, and our ongoing work features this extension. The main challenges will most likely come from a potentially more complex skeleton set structure and from implementation considerations in more complex geometries.

Having such fast direct solvers available is of most relevance to solve large scale boundary value and evolution PDE problems, with applications such as fluid-structure interaction and 3D scattering problems.

\item{\textbf{Parallel implementation:}} as we have noted, the algorithms and implementation presented are serial. However, our direct solver is ideally suited for parallel implementation: intensive linear algebra operations are performed for each node on a binary tree, and the most expensive computations (inverse compression) require only one upward pass.

\item{\textbf{Direct solvers for oscillatory problems:}} ultimately, fast direct solvers such as this are needed to tackle oscillatory problems arising from accoustics and electromagnetics (e.g. Helmholtz equation and Maxwell equations) in moderate and high frequency regimes.
\end{itemize}

\paragraph{Acknowledgments:} We would like to thank Mark Tygert and Leslie Greengard's group for many helpful insights during the design and implementation of our algorithm. 

\bibliography{HSS}

\end{document}